\documentclass[12pt]{article}
\usepackage{amssymb}
\usepackage{graphicx}
\usepackage{caption}
\usepackage{subcaption}

\usepackage{amsfonts}
\usepackage{latexsym}
\usepackage{amsthm}
\usepackage{amsmath}

\usepackage{verbatim}
\usepackage{float}

\usepackage{amssymb, amscd}

\usepackage{amsmath,amsfonts,graphicx,url,amscd}

\usepackage{url}
\usepackage[T1]{fontenc}
\usepackage{color}

  \usepackage{amsmath,amssymb,amsthm,amsfonts}
  \usepackage{pgfplots}
  \usepackage{tikz-3dplot}
  \usepackage{hyperref}
  \usepackage{graphicx}
  \usepackage{tikz}
  \usepackage{mathtools}

\usetikzlibrary{decorations.pathreplacing,decorations.markings,shapes,arrows,arrows.meta,decorations.markings,calc}
\usetikzlibrary{intersections}
\tikzset{->-/.style={
	decoration = {
    markings,
    mark = at position .535 with {\arrow[scale=1.2]{latex}}
    },
    postaction = {decorate}
}}
\tikzset{mid arrow/.style={postaction={decorate,decoration={
        markings,
        mark=at position .5 with {\arrow[#1]{stealth}}
      }}}
}

\newtheorem{problem}{Problem}
\newtheorem{theorem}[problem]{Theorem}
\newtheorem{remark}[problem]{Remark}
\newtheorem{prob}[problem]{Problem}
\newtheorem{defin}[problem]{Definition}
\newtheorem{prop}[problem]{Proposition}

\newtheorem{lema}[problem]{Lemma}
\newtheorem{exam}[problem]{Example}

\usepackage[normalem]{ulem}

\title{{Cyclohedron and Kantorovich-Rubinstein polytopes}}
\author{{Filip D. Jevti\'{c}} \\ {\small Department of Mathematical Sciences}\\[-2mm] {\small University of Texas at Dallas}\\[-2mm] {\small Mathematical Institute SASA, Belgrade}
\and {Marija Jeli\'{c}}\\ {\small Faculty of Mathematics}\\[-2mm] {\small University of Belgrade}
\and Rade T. \v Zivaljevi\' c\\ {\small Mathematical Institute}\\[-2mm] {\small SASA,   Belgrade}\\[-2mm]}

\begin{document}
\date{March 30, 2018}
 \title{{Cyclohedron and Kantorovich-Rubinstein polytopes}}

\author{{Filip D. Jevti\'{c}} \\ {\small Department of Mathematical Sciences}\\[-2mm] {\small University of Texas at Dallas}\\[-2mm] {\small Mathematical Institute SASA}
\and {Marija Jeli\'{c}}\\ {\small Faculty of Mathematics}\\[-2mm] {\small University of Belgrade}
\and Rade T. \v{Z}ivaljevi\'{c}\\ {\small Mathematical Institute}\\[-2mm] {\small SASA,   Belgrade}\\[-2mm]}

\maketitle
\begin{abstract}\noindent
We show that the cyclohedron (Bott-Taubes polytope) $W_n$ arises as the polar dual of a   Kantorovich-Rubinstein polytope $KR(\rho)$, where $\rho$ is an explicitly described  quasi-metric (asymmetric distance function) satisfying strict triangle inequality. From a broader perspective, this phenomenon  illustrates the relationship between a nestohedron $\Delta_{\mathcal{\widehat{F}}}$ (associated to a building set $\mathcal{\widehat{F}}$) and its non-simple deformation $\Delta_{\mathcal{F}}$, where $\mathcal{F}$ is an {\em irredundant} or {\em tight basis} of $\mathcal{\widehat{F}}$ (Definition~\ref{def:irredundant}). Among the consequences are a new proof of a recent result of Gordon and Petrov (\textit{Arnold Math.\ J.}\, \textbf{3} (2), 205--218 (2017)) about $f$-vectors of generic Kantorovich-Rubinstein polytopes and an extension of a theorem of Gelfand, Graev, and Postnikov, about triangulations of the type A, positive root polytopes.
\end{abstract}

\renewcommand{\thefootnote}{}
\footnotetext{This research was supported by the Grants 174020 and
174034 of the Ministry of Education, Science and Technological Development
of Serbia.}

\section{Introduction}
\label{sec:intro}

Motivated by the classic Kantorovich-Rubinstein theorem,  A.M.~Vershik described in \cite{v15} a canonical correspondence between finite metric spaces $(X,\rho)$ and convex polytopes in the vector space  $V_0(X)\subset \mathbb{R}^X$ of all signed measures on $X$ with total mass equal to $0$. More explicitly, each finite metric space $(X,\rho)$ is associated a {\em fundamental polytope} $KR(\rho)$ (Kantorovich-Rubinstein polytope) spanned by $e_{x,y} = \frac{e_x - e_y}{\rho(x,y)}$ where $\{e_x\}_{x\in X}$ is the canonical basis in $\mathbb{R}^X$.

\medskip
Kantorovich-Rubinstein polytope $KR(\rho)$ can be also described as the dual of the {\em  Lipschitz polytope} $Lip(\rho)$ where,
\begin{equation}\label{eqn:Lip}
  Lip(\rho) = \{f \in \mathbb{R}^X \mid (\forall x,y\in X)\, f(x) - f(y)\leqslant \rho(x,y)\},
\end{equation}
and two functions $f, g\in \mathbb{R}^X$ are considered equal if they differ by a constant.

\medskip
Vershik raised in \cite{v15} a general problem of studying (classifying) finite metric spaces according to the combinatorics of their fundamental polytopes.

\smallskip
Gordon and Petrov in a recent paper \cite{gp} proved a very interesting result that the $f$-vector of the Kantorovich-Rubinstein polytope $KR(\rho)$ is one and the same for all  sufficiently generic metrics on $X$. They obtained this result as a byproduct of a careful combinatorial analysis of face posets of Lipschitz polytopes. The invariance of the $f$-vector of $KR(\rho)$ can be also deduced from the fact that the {\em type A root polytope} $Root_n := {\rm Conv}(\mathcal{A}_{n})$, where $\mathcal{A}_{n} = \{e_i-e_j \mid 1\leqslant i\neq j\leqslant n\}$, is unimodular in the sense of  \cite[Definition~6.2.10]{delta} (see also the outline in Section~\ref{sec:unimod}).

\bigskip
Our point of departure was an experimentally observed fact that the {\em generic $f$-vector} of Gordon and Petrov coincides with the $f$-vector of  (the dual of) the cyclohedron (Bott-Taubes polytope) $W_n$. At first sight this is an unexpected phenomenon since $W_n^\circ$ itself is not centrally symmetric and therefore cannot arise as a Kantorovich-Rubinstein polytope $KR(\rho)$  (unless $\rho$ is allowed to be a  quasi-metric!).

\medskip
The symmetry of a metric  is a standard assumption in the usual formulations of the Kantorovich-Rubinstein theorem, see for example \cite[Section~1.2.]{Vil}. However this condition is not necessary. (The proof of this fact is implicit  in \cite[Section~5]{Vil2}, see {\em Particular Case~5.4.} on page 68.) More importantly the `radial vertex perturbation' (Section~\ref{sec:root-prelim}) of a metric may affect its symmetry, so the extension of $KR(\rho)$ to quasi-metrics may be  justified both by the  `optimal transport' and the `convex polytopes' point of view.

\medskip
We prove two closely related results which both provide explanations why the cyclohedron $W_n$ (and its dual polytope $W_n^\circ$) appear in the context of generic Kantorovich-Rubinstein polytopes and triangulations of the type A root polytope.

\medskip
In the first result (Theorem~\ref{thm:glavna-A}), we construct a  map $\phi_n : W_n^\circ \rightarrow Root_n$ which is simplicial on the boundary $\partial(W_n^\circ)$ and maps bijectively $\partial(W_n^\circ)$ to the boundary
 $\partial(Root_n)$ of the root polytope. (In particular we obtain a triangulation of $\partial(Root_n)$ parameterized by faces of $W_n$.)

 This construction is purely combinatorial and diagrammatic in nature. It relies on a combinatorial description of $W_n$ as a graph associahedron \cite{devadoss} and describes simplices in $\partial(W_n^\circ)$ as {\em admissible families} of intervals (arcs) in the cycle graph $C_n$.

 \smallskip
 Theorem~\ref{thm:glavna-A} can also be seen as an extension of a result of Gelfand, Graev, and Postnikov \cite[Theorem~6.3.]{ggp} who described a coherent triangulation of the type A, positive root polytope $Root_n^+ = {\rm Conv}\{e_i-e_j \mid 1\leqslant i\leqslant j\leqslant n\}$. For illustration, the {\em standard trees} depicted in \cite[Figure~6.1.]{ggp}
 may be interpreted as our {\em admissible families of arcs} (as exemplified in Figure~\ref{fig-arc-1}) where all arcs are oriented from left to right.

 \medskip
  In the second result (Theorem~\ref{thm:glavna-B} and Proposition~\ref{prop:quasi-metric-cyclo}) we prove the existence and than explicitly construct a canonical quasi-metric $\rho$ such that the associated Kantorovich-Rubinstein polytope $KR(\rho)$ is a geometric realization of the polytope $W_n^\circ$.

  \smallskip
  This result has a more geometric flavor since it relies on a {\em nestohedron representation} \cite{po,fs} of the cyclohedron as the Minkowski sum $W_n = \Delta_{\mathcal{\widehat{F}}} = \sum_{F\in \mathcal{\widehat{F}}}~\Delta_F$ of simplices. In this approach the relationship between the cyclohedron $W_n$ and the dual $(Root_n)^\circ$ of the root polytope is seen as a special case of a more general construction linking a nestohedron $\Delta_{\mathcal{\widehat{F}}}$ and its Minkowski summand $\Delta_{\mathcal{F}}$, where  $\mathcal{\widehat{F}}$  is a building set and $\mathcal{F}$ its {\em irredundant basis} (Definition~\ref{def:irredundant}).

\medskip
In Section~\ref{sec:alternative} we briefly outline a different plan (suggested by a referee) for constructing quasi-metrics of ``cyclohedral type''. This approach relies on the analysis of the combinatorial structure of Lipschitz polytopes for generic measures, as developed in \cite{gp}.

\medskip

In `Concluding remarks' (Section~\ref{sec:concluding})  we discuss the significance of Theorems~\ref{thm:glavna-A} and \ref{thm:glavna-B}.  For example we demonstrate (Section~\ref{sec:GP-concl}) how the motivating result of Gordon and Petrov \cite[Theorem~1]{gp} can be deduced from the known results about the $f$-vectors of cyclohedra. We also offer a glimpse into potentially interesting future developments including  the study of `tight pairs' $(\mathcal{\widehat{F}}, \mathcal{F})$  of hypergraphs (Section~\ref{sec:tight}) and the `canonical quasitoric manifolds' associated to {\em combinatorial quasitoric pairs} $(W_n, \phi_n)$ (Section~\ref{sec:canonical-quasitoric}).


\section{Preliminaries}

\subsection{Kantorovich-Rubinstein polytopes}

Let $(X,\rho)$, $\vert X\vert =n$, be a finite metric space and let $V(X) := \mathbb{R}^X \cong \mathbb{R}^n$ be the associated vector space of real valued functions (weight distributions, signed measures) on $X$. In particular, $V_0(X) := \{\mu\in V(X) \mid \mu(X)=0\}$ is the vector subspace of measures with total mass equal to zero, while $\Delta_X := \{\mu\in V(X)\mid \mu(X) = 1  \wedge (\forall x\in X)\, \mu(\{x\})\geqslant 0\}$ is the simplex of probability measures.

\medskip
Let $\mathcal{T}_\rho(\mu, \nu)$ be the cost of optimal transportation of measure
$\mu$ to measure $\nu$, where the cost of transporting the unit mass from $x$ to $y$ is $\rho(x,y)$. Then \cite{v13, Vil}, there exists a norm $\|\cdot \|_{KR}$ on $V_0(X)$ (called the Kantorovich-Rubinstein norm), such that,
\[
             \mathcal{T}_\rho(\mu, \nu) = \|\mu - \nu \|_{KR},
\]
for each pair of probability measures $\mu, \nu\in \Delta_X$. By definition, the Kantorovich-Rubinstein polytope $KR(\rho)$, or the {\em fundamental polytope} \cite{v15}, associated to $(X,\rho)$, is the corresponding unit ball in $V_0(X)$,
\begin{equation}\label{eqn:KR-poly}
  KR(\rho) = \{x\in V_0(X) \mid \|x \|_{KR} \leqslant 1\}.
\end{equation}
The following explicit description for $KR(\rho)$ can be deduced from the Kantorovich-Rubinstein theorem (Theorem~1.14 in \cite{Vil}),
\begin{equation}\label{eqn:KR-conv}
  KR(\rho) =  {\rm Conv}\left\{  \frac{e_x - e_y}{\rho(x,y)} \mid x\in X \right\},
\end{equation}
where $\{e_x\}_{x\in X}$ is the canonical basis in $\mathbb{R}^X$.

\medskip\noindent
\begin{prob}{\rm (A.M.~Vershik \cite{v15})}
  Study and classify metric spaces according to combinatorial properties of their Kantorovich-Rubinstein polytopes.
\end{prob}

\subsection{Root polytopes}\label{sec:root-prelim}

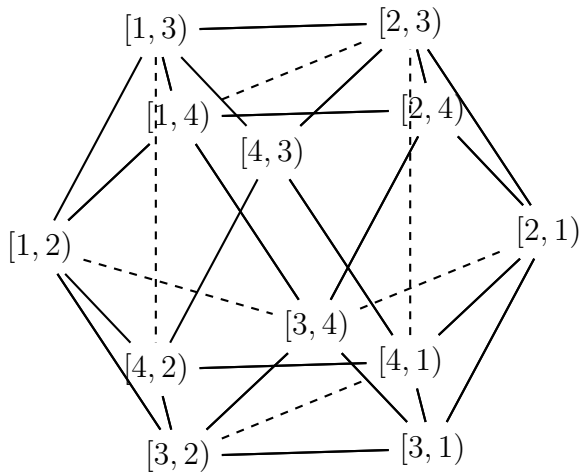
\begin{figure}[htb]
	\centering
	\tdplotsetmaincoords{110}{-40}
    \begin{tikzpicture}[tdplot_main_coords,scale=0.6]

       \coordinate (e0) at (-4,4,4);
       \coordinate (e1) at (4,-4,4);
       \coordinate (e2) at (4,4,-4);
       \coordinate (e3) at (-4,-4,-4);
       \coordinate (f0) at (4,-4,-4);
       \coordinate (f1) at (-4,4,-4);
       \coordinate (f2) at (-4,-4,4);
       \coordinate (f3) at (4,4,4);

       \node (a03) at (0,4,4) {$\left[1,4\right)$};
       \node (a23) at (4,4,0) {$\left[3,4\right)$};
       \node (a21) at (0,4,-4) {$\left[3,2\right)$};
       \node (a01) at (-4,4,0) {$\left[1,2\right)$};
       \node (a32) at (-4,-4,0) {$\left[4,3\right)$};
       \node (a12) at (0,-4,4) {$\left[2,3\right)$};
       \node (a30) at (0,-4,-4) {$\left[4,1\right)$};
       \node (a10) at (4,-4,0) {$\left[2,1\right)$};
       \node (a13) at (4,0,4) {$\left[2,4\right)$};
       \node (a20) at (4,0,-4) {$\left[3,1\right)$};
       \node (a31) at (-4,0,-4) {$\left[4,2\right)$};
       \node (a02) at (-4,0,4) {$\left[1,3\right)$};

		
       	\draw[color=black,thick] (a30) -- (a31) -- (a21) -- (a20) --(a30);
        \draw[color=black,thick] (a30) -- (a10) -- (a12) -- (a32) -- (a30);
        \draw[color=black,thick] (a02) -- (a32) -- (a12) -- (a02);
        \draw[color=black,thick] (a12) -- (a13) -- (a10) -- (a12);
        \draw[color=black,thick] (a32) -- (a30) -- (a31) -- (a32);
        \draw[color=black,thick] (a02) -- (a03) -- (a01) -- (a02);
        \draw[color=black,thick] (a03) -- (a13) -- (a23) -- (a03);
        \draw[color=black,thick] (a30) -- (a20) -- (a10) -- (a30);
        \draw[color=black,thick] (a23) -- (a20) -- (a21) -- (a23);
        \draw[color=black,thick] (a01) -- (a31) -- (a21) -- (a01);
        \draw[color=black,thick] (a02) -- (a03) -- (a13) -- (a12) -- (a02);
        \draw[color=black,thick] (a13) -- (a10) -- (a20) -- (a23) -- (a13);
        \draw[color=black,thick] (a01) -- (a21) -- (a23) -- (a03) -- (a01);

        \draw[color=black,thick,dashed] (a02) -- (a31);
        \draw[color=black,thick,dashed] (a01) -- (a23);
        \draw[color=black,thick,dashed] (a21) -- (a30);
        \draw[color=black,thick,dashed] (a23) -- (a10);
        \draw[color=black,thick,dashed] (a03) -- (a12);
        \draw[color=black,thick,dashed] (a30) -- (a12);

   \end{tikzpicture} 
 \caption{The boundary triangulation of $Root_4$ described in Theorem~\ref{thm:glavna-A}.}
\label{fig-Root4-1}
\end{figure}

The convex hull of the roots of a classical root system is called a root polytope.
In particular the {\em type A root polytope}, associated to the root system of type $A_{n-1}$, is the following polytope (Fig.\ref{fig-Root4-1}),
\begin{equation}\label{eqn:root-poly-typeA}
  Root_n = {\rm Conv}\{e_i-e_j \mid 1\leqslant i \neq j\leqslant n\}.
\end{equation}
It immediately follows from (\ref{eqn:root-poly-typeA}) that the root polytope admits the following Minkowski sum decomposition,
\begin{equation}\label{eqn:root-Mink}
  Root_n = \Delta + \nabla = \Delta + (-\Delta) = \Delta - \Delta,
\end{equation}
where $\Delta = \Delta_n = {\rm Conv}\{e_i\}_{i=1}^n$ and $\nabla = -\Delta = {\rm Conv}\{-e_i\}_{i=1}^n$.

\medskip
By definition, $Root_n$ is the Kantorovich-Rubinstein polytope associated to the metric $\rho$ where $\rho(x,y) =1$ for each $x\neq y$. Conversely, in light of (\ref{eqn:KR-conv}), each Kantorovich-Rubinstein polytope can be seen as a {\em radial, vertex perturbation} of the root polytope $Root_n$.

\subsection{Unimodular triangulations and equidecomposable polytopes}
\label{sec:unimod}

 A triangulation of a convex polytope $Q$ is tacitly assumed to be without new vertices. A triangulation of the boundary sphere $\partial(Q)$ of $Q$ is referred to as a {\em boundary triangulation}. Each triangulation of $Q$ produces the associated boundary triangulation (but not the other way around).

  \smallskip
  The $f$-vector of a triangulation is the $f$-vector of the associated simplicial complex. Different triangulations of either the polytope $Q$ or its boundary $\partial(Q)$ may have different face numbers, so in general the $f$-vector of a triangulation is not uniquely determined by the polytope $Q$. The simplest examples illustrating this phenomenon are the bipiramid over a triangle and the $3$-dimensional cube (the latter admits triangulations with both $5$ and $6$, three dimensional simplices).

 \smallskip
 The polytopes, all of whose triangulations have the same face numbers ($f$-vectors), are called {\em equidecomposable}, see \cite{ba} or  \cite[Section~8.5.3]{delta}. A notable class of equidecomposable polytopes are lattice polytopes which are {\em unimodular} in the sense that each full dimensional simplex spanned by its vertices has the same volume, see Definition~6.2.10  and Section~9.3 in \cite{delta}. Unimodularity of a polytope  immediately implies that the top dimensional face numbers are independent of a triangulation. In light of Theorem~8.5.19. from \cite{delta}, this condition guarantees that the polytope is equidecomposable, i.e.\ that the $f$-vector is the same for all triangulations.

 \medskip
  A notable example of an equidecomposable polytope is the product of two simplices, see \cite[Section~6.2]{delta}. As a consequence of (\ref{eqn:root-Mink}), each face of the root polytope $Root_n$ is a product of two simplices. From here we immediately deduce that all boundary triangulations of $Root_n$ have the same $f$-vector.

  \medskip
  Gordon and Petrov \cite{gp} observed that each Kantorovich-Rubinstein polytope $KR(\rho)$, for a sufficiently generic metric $\rho$, induces a {\em regular boundary triangulation} of the root polytope $Root_n$. This observation allowed them to determine the $f$-vector of a generic K-R polytope, and to obtain some other qualitative and quantitative information about these polytopes.

  \medskip
  Our Theorem~\ref{thm:glavna-B} identifies this $f$-vector as the $f$-vector of the polytope $W_n^\circ$, dual to the $f$-vector of an $(n-1)$-dimensional cyclohedron.

\subsection{Kantorovich-Rubinstein polytopes for quasi-metrics}

  Each Kantorovich-Rubinstein polytope associated to a metric $\rho$ is centrally symmetric (as a consequence of the symmetry $\rho(x,y) = \rho(y,x)$ of the metric $\rho$).
  The cyclohedron $W_n$ is not centrally symmetric, so it is certainly not one of the Kantorovich-Rubinstein polytopes.  However, as a consequence of Theorem~\ref{thm:glavna-B}, it arises as a {\em generalized K-R polytope} associated to a not necessarily symmetric distance function (quasi-metric).

  \begin{defin}
    A non-negative function $\rho : X\times X\rightarrow \mathbb{R}^+$ is a {\em quasi-metric} (asymmetric distance function) if,
    \begin{enumerate}
      \item    $(\forall  x,y\in X)\, (\rho(x,y) =0  \, \Leftrightarrow \, x=y)$;
      \item   $(\forall x,y,z\in X) \, \rho(x,z) \leqslant \rho(x,y) + \rho(y,z)$.
      \end{enumerate}
     The Kantorovich-Rubinstein polytope $KR(\rho)$, associated to a quasi-metric $\rho$, is defined by the same formula (\ref{eqn:KR-conv}) as its symmetric counterpart.
  \end{defin}
Many basic facts remain true for generalized K-R polytopes. For illustration, here is a result which extends (with the same proof) a result of  Melleray et al. \cite[Lemma~1]{mpv}.

\begin{prop}
Let $X$ be a finite set. Assume that $\rho:X\times X \rightarrow \mathbb{R}_{\geqslant 0}$ is a non-negative function such that $\rho(x,y) = 0$ if and only if $x=y$. Let $KR(\rho)$ be the polytope defined by the equation (\ref{eqn:KR-conv}). Then $\rho$ is a quasi-metric on $X$ if and only if none of the points $e_{x,y} = \frac{e_x-e_y}{\rho(x,y)}$  (for $x\neq y$)  is in the interior of $KR(\rho)$.
\end{prop}

\section{Preliminaries on the cyclohedron $W_n$}\label{sec:cyclo-prelim}

\subsection{Face lattice of the cyclohedron $W_n$ }
\label{sec:face}

The face lattice $\mathcal{F}(W_n)$ of the $(n-1)$-dimensional {\em cyclohedron} $W_n$  (Bott-Taubes polytope) admits two closely related combinatorial description.

\medskip
In the first description  \cite{sta97}, similar to the description of the $(n-1)$-dimensional associahedron $K_n$ (Stasheff polytope),  the lattice $\mathcal{F}(W_n)$ arises as the collection of all partial cyclic bracketing of a word  $x_1x_2\dots x_n$.

\medskip
Carr and Devadoss \cite{carr-devadoss-1}, in a more general approach, view both polytopes $K_n$ and $W_n$ as instances of the so called {\em graph associahedra}. In this approach, $W_n$ is described as the graph associahedron $\mathcal{P}\Gamma$ corresponding  to the graph $\Gamma = C_n$ (cycle on $n$ vertices), where $\mathcal{F}(W_n)$ is the collection of all {\em valid tubings} on $\Gamma$, see \cite{carr-devadoss-1, devadoss}  for details.

\medskip
The equivalence of the `bracketing' and `tubing' description is easily established, see for example \cite[Lemma~2.3]{carr-devadoss-1} or \cite[Lemma~1.4]{mar99}.
Recall that `graph associahedra' are a specialization of {\em nestohedra}, see Feichtner-Sturmfels \cite{fs}, Postnikov \cite{po}, or Buchstaber-Panov \cite[Section~1.5.]{bp}. In this more general setting, the `valid tubings' appear under the name of `nested sets' associated to a chosen `building set'. A related class of polytopes was studied from a somewhat different perspective by Do\v{s}en and Petri\'{c} in \cite{dp}.

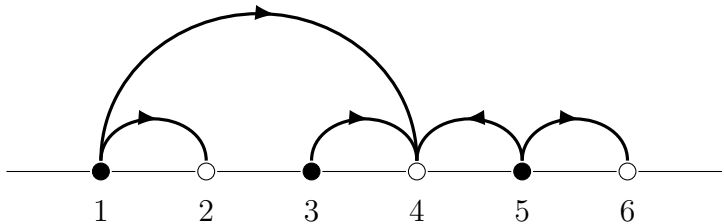
\begin{figure}[htb]
	\centering
	\begin{tikzpicture}[scale=1.4]
        \node (l) at (-1,0) {};
        \node (r) at (6,0) {};
      	\node (p0) at (0,0) {};
	    \node (p1) at (1,0) {};
		\node (p2) at (2,0) {};
      	\node (p3) at (3,0) {};
    	\node (p4) at (4,0) {};
    	\node (p5) at (5,0) {};
		\draw[fill=black] (p0) circle(0.08) node[yshift=-15pt] {$1$};
        \draw (p1) circle(0.08) node[yshift=-15pt] {$2$};
        \draw[fill=black] (p2) circle(0.08) node[yshift=-15pt] {$3$};
        \draw (p3) circle(0.08) node[yshift=-15pt] {$4$};
        \draw[fill=black] (p4) circle(0.08) node[yshift=-15pt] {$5$};
        \draw (p5) circle(0.08) node[yshift=-15pt] {$6$};

        \draw (l) -- (p0) -- (p1) -- (p2) -- (p3) -- (p4) -- (p5) -- (r);

      	\draw[very thick,->-] (p0) to[out=90,in=180] (1.5,1.5) to[out=0,in=90] (p3);
        \draw[very thick,->-] (p0) to[out=90,in=180] (.5,.5) to[out=0,in=90] (p1);
        \draw[very thick,->-] (p2) to[out=90,in=180] (2.5,.5) to[out=0,in=90] (p3);
        \draw[very thick,->-] (p4) to[out=90,in=0] (3.5,.5) to[out=180,in=90] (p3);
        \draw[very thick,->-] (p4) to[out=90,in=180] (4.5,.5) to[out=0,in=90] (p5);
     \end{tikzpicture}
\caption{Admissible family of oriented arcs.}
\label{fig-arc-2}
\end{figure}

\bigskip
We use in this paper a slightly modified description of the poset $\mathcal{F}(W_n)$ which allows us to use pictorial description of `valid tubings' (partial bracketings), see Fig.\ref{fig-arc-2} for an example. A similar description was used by Gelfand, Graev, and Postnikov \cite{ggp}, where these pictorial representations appeared in the form of the so called `standard trees', see \cite[Section~6]{ggp}.

 The vertex set of the cycle graph $C_{n}$   is the set $V(C_{n}) =  [n] := \{1, 2, \dots, n-1, n\}$ of vertices of a regular $n$-gon, inscribed in a unit circle $S^1$.  We adopt a (counterclockwise) circular order $\prec$ (respectively $\preceq$) on the circle $S^1$,
so in particular $[x,y] = \{z\in S^1\mid x\preceq z \preceq y\}$ is a closed arc (interval) in $S^1$ (similarly $[x,y) = \{z\in S^1\mid x\preceq z \prec y\}$, etc.).
(By convention, $[x,x]=\{x\}$ and $[x,x)=\emptyset$.)

 \medskip By definition  $[x,y]^0 := [x,y]\cap [n]$  (similarly, $[x,y)^0 := [x, y)\cap [n]$) are intervals restricted to the set $[n]$ of vertices of the $n$-gon. If $i\neq j$  are two distinct vertices (elements of $V(C_n) = [n]$ ), then $[i,j-1]^0$ is precisely the tube (in the sense of \cite{carr-devadoss-1}) associated to the interval $[i,j)$. Conversely, each tube $[u,v]^0$ (a proper connected component in the graph $C_n$) is associated a half-open interval $[u,v+1)$ in the circle $S^1$. A moment's reflection reveals that each valid tubing (in the sense of \cite{carr-devadoss-1}) corresponds to an admissible family  of half-open intervals, in the sense of the following definition.

\begin{defin}\label{def:admissible}
A collection $\mathcal{T} = \{I_1, I_2, \dots, I_k\}$ of half-open intervals $I_j = [a_j, b_j)$ (where $\{a_j, b_j\}\subset [n]$ and $a_j\neq b_j$ for each $j$) is {\em admissible} if for each  $1\leqslant i < j \leqslant k$, exactly one of the following two alternatives is true,
 \begin{enumerate}
   \item   If $I_i\cap I_j \neq \emptyset$ then $I_i, I_j$ are comparable in the sense that either $I_i\subsetneq I_j$ or $I_j\subsetneq I_i$;
   \item  $I_i\cap I_j = \emptyset$ and $I_i\cup I_j$ is not an interval (meaning that neither $b_j=a_i$ nor $b_i=a_j$).
 \end{enumerate}
 \end{defin}

 \begin{prop}\label{prop:admissible}
   The face lattice $\mathcal{F}(W_n)$ of the $(n-1)$-dimensional {\em cyclohedron} $W_n$ is isomorphic to the poset of all  admissible collections  $\mathcal{T} = \{I_j\}_{j=1}^k$ of half-open intervals in $S^1$ with endpoints in $[n]$. Individual arcs (half-open intervals) correspond to facets of $W_n$ while the empty set is associated to the polytope $W_n$ itself.
 \end{prop}

\begin{remark}\label{rem:admissible}{\rm
 The dual $W_n^\circ$ of the cyclohedron is a simplicial polytope. It follows from Proposition~\ref{prop:admissible} that the face poset of the boundary  $\partial(W_n^\circ)$ of $W_n^\circ$ is isomorphic to the simplicial complex with vertices $V = \{[i,j)\mid 1\leqslant i\neq j\leqslant n\}$ (all half-open intervals with endpoints in $[n]$) where $\mathcal{T}\subset V$ is a simplex if and only if $\mathcal{T}$ is an admissible family of half-open intervals (Definition~\ref{def:admissible}).
}
\end{remark}

\begin{figure}[h]
	\centering
	\begin{subfigure}[b]{0.4\textwidth}
    	\centering
		\begin{tikzpicture}[scale=1.4]
  	\def \N {4} 
    \def \s {1.1} 
    \def \E {2/1,3/1,3/0}

	\pgfmathsetmacro{\tmp}{\N-1}
    \pgfmathsetmacro{\BEG}{0};
    \pgfmathsetmacro{\END}{\tmp*\s};

	\draw (\BEG,0) -- (\END,0);
    \foreach \x in {0,...,\tmp}{
       	\node[draw,circle,fill=white,inner sep=0pt, minimum size=5pt] (n\x) at ({\x*\s},0) {};
    }
    \foreach \x in {1,...,\N}{
    	\node (l\x) at ({(\x-1)*\s},-.25) {$\x$};
    }
    \foreach \i/\j in \E {
        \node[draw,circle,fill=black,inner sep=0pt, minimum size=5pt] (n\i) at ({\i*\s},0) {};
        \pgfmathsetmacro{\XC}{((\i+\j)/2)*\s};
        \pgfmathsetmacro{\YC}{(abs(\j-\i)/2)*\s};

        \pgfmathparse{(\j - \i) > 0.001 ? int(1) : int(0)};
    	\ifnum\pgfmathresult=1
			\draw[very thick,->-] (n\i) to[out=90,in=180] (\XC,\YC) to[out=0,in=90] (n\j);
        \else
            \draw[very thick,->-] (n\i) to[out=90,in=0] (\XC,\YC) to[out=180,in=90] (n\j);
        \fi
	}
\end{tikzpicture}
        \caption{$\mathcal{T}_1 = \{[3,2), [4,2), [4,1)\}$}
	\end{subfigure}
    \begin{subfigure}[b]{0.4\textwidth}
    	\centering
		\begin{tikzpicture}[scale=1.4]
	\def \N {4} 
    \def \s {1.1} 
    \def \E {2/1,2/0,3/0}

    \pgfmathsetmacro{\tmp}{\N-1}
  	\pgfmathsetmacro{\BEG}{0};
    \pgfmathsetmacro{\END}{\tmp*\s};

	\draw (\BEG,0) -- (\END,0);
    \foreach \x in {0,...,\tmp}{
    	\node[draw,circle,fill=white,inner sep=0pt, minimum size=5pt] (n\x) at ({\x*\s},0) {};
    }
    \foreach \x in {1,...,\N}{
        \node (l\x) at ({(\x-1)*\s},-.25) {$\x$};
    }
    \foreach \i/\j in \E {
        \node[draw,circle,fill=black,inner sep=0pt, minimum size=5pt] (n\i) at ({\i*\s},0) {};
        \pgfmathsetmacro{\XC}{((\i+\j)/2)*\s};
        \pgfmathsetmacro{\YC}{(abs(\j-\i)/2)*\s};

		\pgfmathparse{(\j - \i) > 0.001 ? int(1) : int(0)};
    	\ifnum\pgfmathresult=1
			\draw[very thick,->-] (n\i) to[out=90,in=180] (\XC,\YC) to[out=0,in=90] (n\j);
  		\else
            \draw[very thick,->-] (n\i) to[out=90,in=0] (\XC,\YC) to[out=180,in=90] (n\j);
   		\fi
     }
\end{tikzpicture}
        \caption{$\mathcal{T}_2 = \{[3,2), [3,1), [4,1)\}$}
	\end{subfigure}
\caption{Two facets of $W_4^\circ$ (vertices of $W_4$) .}
\label{fig-arc-1}
\end{figure}
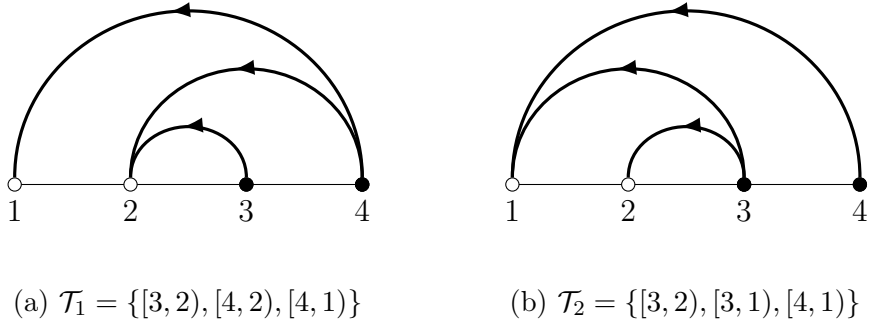

The following proposition shows that admissible families (in the sense of Definition~\ref{def:admissible}) can be naturally interpreted as directed trees (directed forests)

\begin{prop}\label{prop:tree}
  Each admissible family $\mathcal{T} = \{[i_\nu, j_\nu)\}_{\nu=1}^k$ of arcs can be  interpreted as a digraph, where $[i_\nu, j_\nu)\in \mathcal{T}$ defines an oriented edge from $i_\nu$ to $j_\nu$ (as in Figures~\ref{fig-arc-2} and \ref{fig-arc-1}). If the orientation of arcs is neglected, than we obtain an ordinary graph $\Gamma_\mathcal{T}$. We claim that for each admissible family $\mathcal{T}$, the associated graph $\Gamma_\mathcal{T}$ has no cycles.
\end{prop}

\medskip\noindent
{\bf Proof:}
Indeed, suppose that $[u_1, v_1), (v_1, u_2], [u_2, v_2), \dots, [u_s, v_s), (v_s,u_1] $
is a minimal cycle in $\Gamma_\mathcal{T}$. We may assume without loss of generality that $u_1 \prec v_1\prec u_2$ (in the counterclockwise circular order  on $S^1$). From here we deduce that remaining indices also follow the circular order,
\begin{equation}\label{eqn:circ}
   u_1\prec v_1 \prec u_2 \prec v_2 \prec \dots \prec u_s \prec v_s \prec u_1,
\end{equation}
otherwise two different arcs would cross (which would violate the assumption that $\mathcal{T}$ is admissible). Moreover, for the same reason, the sequence (\ref{eqn:circ}) winds around the circle $S^1$ only once.
This however leads to a contradiction since the intervals $I = (v_1, u_2]$ and $J = (v_2, u_3]$ would have a non-empty intersection, while neither $I\subset J$ nor $J\subset I$ (a contradiction with Definition~\ref{def:admissible}). \hfill $\square$

\begin{defin}\label{def:circle_interval}
Let $a=\left[ i,j \right)$ be a half-open circle interval. By definition $s(a)=i$ is the source of $a$ and $t(a)=j$ is the sink or the terminal point of $a$. For an admissible family
 $\mathcal{T}=\left\lbrace \left[ i_\nu, j_\nu\right) \mid \nu\in \left[ k \right]  \right\rbrace$ of intervals (arcs), the associated source and sink sets are,
\[
s(\mathcal{T})=\left\lbrace i_\nu \mid \nu\in \left[ k \right]  \right\rbrace  \quad \mbox{ {\rm and} } \quad
t(\mathcal{T})=\left\lbrace j_\nu \mid \nu\in \left[ k \right]  \right\rbrace.
\]
\end{defin}

Note that, as a consequence of Definition~\ref{def:admissible}, $s(\mathcal{T})\cap t(\mathcal{T})=\emptyset$ for each admissible family $\mathcal{T}$ of intervals.

\subsection{Automorphism group of the cyclohedron}
\label{sec:auto-group}

Each automorphism of a graph $\Gamma$ induces an automorphism of the associated graph associahedron $\mathcal{P}_\Gamma$.  The group of all automorphisms of the cycle graph $C_n$ is the dihedral group $D_{2n}$ of order $2n$. It immediately follows that both the $(n-1)$-dimensional cyclohedron $W_{n}$ and its polar polytope $W_n^\circ$  are invariant under the action of the dihedral group $D_{2n}$.

\medskip
Let $D_{2n} \cong \langle  \alpha,\beta \vert\,  \alpha^n = \beta^2 = e, \beta\alpha\beta = \alpha^{n-1} \rangle$ be a standard presentation of the group $D_{2n}$ where $\alpha$ is the cyclic permutation of $C_n$ (corresponding to the rotation of the regular polygon through the angle $2\pi/n$) and $\beta$ is the involution (reflection) which keeps the vertex $n$ fixed.

\medskip
Then the action of $D_{2n}$ on $W_n$ and $W_n^\circ$ can be described as follows.

\begin{prop}\label{prop:auto}
  Suppose that  $C_n$ is the cycle graph, realized as a regular polygon inscribed in the unit circle. Let $[i,j)$ be a half-open interval representing a  vertex (face) of the polytope $W_n^\circ$ (respectively polytope $W_n$).  Then $\alpha([i, j)) := [i+1, j+1)$  and $\beta([i, j)) := [n-j+1, n-i+1)$.
\end{prop}

\subsection{Canonical map $\phi_n$}

The associahedron $K_n$ may be described as the {\em secondary polytope} \cite{gkz}, associated to all subdivisions of a convex $(n+2)$-gon by configurations of  non-crossing diagonals. It was shown by R.~Simion  \cite{si} that a similar description exists for $W_n$, provided we deal only with centrally symmetric configurations. The polytopes $K_n$ and $W_n$ are sometimes referred to as the {\em type A} and {\em type B} associahedra. This classification emphasizes a connection with type A or B root systems, the corresponding hyperplane arrangements etc. In this section we relate $W_n$ to the root system of type $A_{n-1}$, in other words $W_n$ may also be interpreted as a `type A associahedron'.

\medskip
Let $\{e_i\}_{i=1}^n$ be the standard basis in $\mathbb{R}^n$ and let $\mathcal{A}_n = \{e_i-e_j\}_{1\leqslant i\neq j\leqslant n}$ be the associated root system of type $A_{n-1}$. The type A root polytope is introduced in Section~\ref{sec:root-prelim} as the convex hull $Root_n = \mbox{\rm conv} \{e_i-e_j \mid 1\leq i\neq j\leq n\}$ of the set $\mathcal{A}_n$ of all roots. (We warn the reader that this terminology may not be uniform, for example the root polytopes introduced in \cite{ggp, po} deal only with the set  $\mathcal{A}_n^+ = \{e_i-e_j\}_{1\leqslant i < j\leqslant n}$  of positive roots.)

\medskip
The following definition introduces a canonical map which links the (dual of the) cyclohedron $W_n$ to the root system of type $A_{n-1}$, via the root polytope $Root_n$. Recall (Proposition~\ref{prop:admissible}) that the boundary $\partial(W_n^\circ)$ of the polytope dual to the cyclohedron is the simplicial complex of all admissible half-open intervals in $S^1$ with endpoints in $[n]$.

\begin{defin}\label{def:canonical-map}
The map,
\begin{equation}\label{eqn:canonical-map}
 \phi_n : \partial(W_n^\circ) \longrightarrow \partial(Root_n)
\end{equation}
is defined as the simplicial (affine) extension of the map which sends the interval $[i,j)$ (vertex of $W_n^\circ$) to $\phi_n([i,j)):= e_i-e_j\in \mathbb{R}^n$.
\end{defin}

\begin{prop}\label{prop:loc-injective}
The map $\phi_n : \partial(W_n^\circ) \rightarrow \partial(Root_n)$, introduced in Definition~\ref{def:canonical-map},  is one-to-one on faces, i.e.\ it sends simplices of $\partial(W_n^\circ)$ to non-degenerate simplices in the boundary $\partial(Root_n)$ of the root polytope.
\end{prop}

\medskip\noindent
{\bf Proof:} A face of $\partial(W_n^\circ)$ corresponds to an admissible family $\mathcal{T} = \{[i_\nu, j_\nu)\}_{\nu=1}^k$. The associated digraph (also denoted by $\mathcal{T}$) is a directed forest (by Proposition~\ref{prop:tree}).
  The associated unoriented graph $\Gamma_\mathcal{T}$ is a bipartite graph with the shores $P = \{i_1,\dots, i_k\}$ and $Q= \{j_1,\dots, j_k\}$ which has no cycles, i.e.\ $\Gamma_\mathcal{T}$ is a forest. The elements of the corresponding set  $\phi_n(\mathcal{T}) = \{e_{i_\nu} - e_{j_\nu}\}_{j=1}^k$ of vectors may be interpreted as some of the vertices of the product of simplices $\Delta_P\times \Delta_Q$. By Lemma~6.2.8 from \cite[Section~6.2]{delta} these vertices are affinely independent. This implies that $\phi_n$ must be one-to-one on $\mathcal{T}$.
  \hfill $\square$

\medskip
Proposition~\ref{prop:loc-injective} is a very special case of Proposition~\ref{prop:glob-injective}, which claims global injectivity of the canonical map $\phi_n$. The following example illustrates one of the main reasons why the $\phi_n$-images of different simplices have disjoint interiors.

\begin{exam}\label{exam:3} {\rm
By inspection of Figure~\ref{fig-Root4-1} (which illustrates the case $n=4$ of Proposition~\ref{prop:glob-injective}), we observe that the images of different triangles (admissible triples) $\mathcal{T}_1$ and $\mathcal{T}_2$, have disjoint interiors. For example let (Fig.~\ref{fig-arc-1}),
\begin{align*}
\mathcal{T}_1 = \{[3,2), [4,2), [4,1)\}  \quad \mbox{\rm and}\quad \mathcal{T}_2 = \{[3,2), [3,1), [4,1)\}.
\end{align*}
Suppose that the interiors of their images have a non-empty intersection, i.e.\ assume  that there is a solution of the equation,
\begin{align*}\label{eqn:3-illustration}
  \alpha_{3,2}(e_3-e_2) + \alpha_{4,2}(e_4-e_2) + \alpha_{4,1}(e_4-e_1) =
  \beta_{3,2}(e_3-e_2)  +  \beta_{3,1}(e_3-e_1)  +  \beta_{4,1}(e_4-e_1)
\end{align*}
By rearranging the terms we obtain,
\begin{align*}
\alpha'(e_4-e_3) + \alpha''(e_3-e_2) + \alpha'''(e_2-e_1) = \beta'(e_4-e_3) + \beta''(e_3-e_2) + \beta'''(e_2-e_1)
\end{align*}
However, this is impossible since $\alpha$'s and $\beta$'s are positive and $\alpha''>\alpha'> \alpha'''>0$ while $\beta''>\beta'''>\beta'>0$. }
\end{exam}

\begin{remark}\label{rem:korisno}{\rm
  The argument used in the previous example is sufficiently general to cover the case of triangles (admissible triples) $\mathcal{T}_1$ and $\mathcal{T}_2$, which share a common edge (as in Figure~\ref{fig-circ-2}). Indeed, by setting $i=1, j=2$ and $n=4$ this case is reduced to Example~\ref{exam:3}.}
\end{remark}

\section{Cyclohedron and the root polytope I}
\label{sec:cyclo-root-I}

The following theorem is together with Theorem~\ref{thm:glavna-B} one of the central results of the paper.  Informally speaking, it says that there exists a triangulation of the boundary of the $(n-1)$-dimensional type A root polytope $Root_n$ parameterized by proper faces of the $(n-1)$-dimensional cyclohedron.

\begin{theorem}\label{thm:glavna-A}
The map $\phi_n : \partial(W_n^\circ)\rightarrow \partial(Root_n)$, introduced in Definition~\ref{def:canonical-map}, is a piecewise linear homeomorphism of boundary spheres of polytopes $W_n^\circ$ and $Root_n$. The map $\phi_n$ sends bijectively vertices of $\partial(W_n^\circ)$ to vertices of the polytope $Root_n$, while higher dimensional faces of $Root_n$ are triangulated by images of simplices from $\partial(W_n^\circ)$.
\end{theorem}

The proof of Theorem~\ref{thm:glavna-A} is given in the following two sections. Its main part is the proof of the injectivity of the canonical map $\phi_n$.

\subsection{Injectivity of the map $\phi_n$ }
\label{sec:injectivity}

Proposition~\ref{prop:loc-injective} can be interpreted as a result claiming  local injectivity of the map $\phi_n : \partial(W_n^\circ)\rightarrow \partial(Root_n)$. Our central result in this section is Proposition~\ref{prop:glob-injective}, which  establishes global injectivity of this map and provides a key step in the proof of Theorem~\ref{thm:glavna-A}.

\medskip
Recall that $\phi_n : \partial(W_n^\circ)\rightarrow \partial(Root_n)$ is defined as the simplicial map such that $\phi_n\left( \left[i,j\right) \right)=e_i-e_j$ for each pair $i\neq j$. More explicitly, if $x \in \partial W_n^\circ$ is a convex combination of arcs (intervals),
\begin{equation}\label{eqn:conv-comb-1}
x=\sum_{ \left[i,j\right) \in \mathcal{T}} \lambda_{i,j}  \left[i,j\right)
\end{equation}
(where $\mathcal{T}$ is the associated admissible family)  then,
\begin{equation}\label{eqn:conv-comb-2}
\phi_n(x)=\sum_{ \left[i,j\right) \in \mathcal{T}} \lambda_{i,j}  (e_i - e_j).
\end{equation}
We will usually assume that the representation (\ref{eqn:conv-comb-1}) is minimal ($x\in int(\mathcal{T})$) which means that the weights $\{\lambda_{i,j}\}$ satisfy the conditions $\sum_{ \left[i,j\right) \in T} \lambda_{i,j}=1$ and $\left( \forall i,j \right) \lambda_{i,j} > 0$.

\begin{figure}[htb]
	\centering
	\begin{tikzpicture}[scale=0.7]
    	\def \N {12} 
        \def \R {3} 
        \def \E {10/1/1,10/5/0,9/1/0,9/5/1}

        \pgfmathsetmacro{\ang}{360/\N}	
        \pgfmathsetmacro{\tmp}{\N-1};

		\draw (0,0) circle(\R);

        \foreach \x in {0,...,\tmp} {
           	\node[draw,circle,fill=white,inner sep=0pt, minimum size=6pt] (n\x) at ({\x*\ang}:{\R}) {};
      	}
        \node () at ({10*\ang}:{\R+.5}) {$n$};
        \node () at ({9*\ang}:{\R+.5}) {$n-1$};
        \node () at ({5*\ang}:{\R+.5}) {$j$};
        \node () at ({1*\ang}:{\R+.5}) {$i$};

       \foreach \i/\j/\t in \E {
        	\node[draw,circle,fill=black,inner sep=0pt, minimum size=6pt] (n\i) at ({\i*\ang}:{\R}) {};
            \pgfmathparse{abs(\j - \i) < (\N/2) ? int(1) : int(0)};
    		\ifnum\pgfmathresult=1
            	\pgfmathparse{\t > 0.5 ? int(1) : int(0)};
                \ifnum\pgfmathresult=1
					\draw[very thick,->-] (n\i) to[bend right] (n\j);
                \else
					\draw[very thick,dashed,->-] (n\i) to[bend right] (n\j);
                \fi
	   		\else
            	\pgfmathparse{\t > 0.5 ? int(1) : int(0)};
                \ifnum\pgfmathresult=1
					\draw[very thick,->-] (n\i) to[bend left] (n\j);
                \else
					\draw[very thick,dashed,->-] (n\i) to[bend left] (n\j);
                \fi
    		\fi
        }
	\end{tikzpicture}
    \caption{Admissible triangles with a common edge.}
    \label{fig-circ-2}
\end{figure}
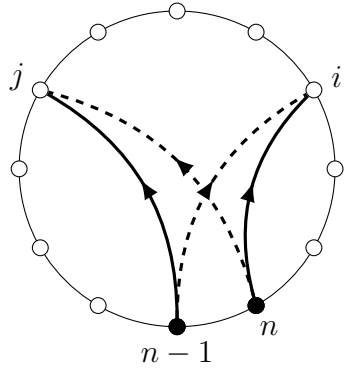

\begin{prop}\label{prop:glob-injective}
The map $\phi_n : \partial(W_n^\circ)\rightarrow \partial(Root_n)$ is injective.
\end{prop}

\medskip\noindent
{\bf Proof:}
Suppose that $\mathcal{T}_1$ and $\mathcal{T}_2$ are two admissible families of intervals (representing two faces of $\partial W_n^{\circ}$). We want to show that $\mathcal{T}_1$ and $\mathcal{T}_2$ must be equal if,
\begin{align}\label{eqn:cond-inj}
\phi_n(int(\mathcal{T}_1)) \cap  \phi_n(int(\mathcal{T}_2)) \neq \emptyset .
\end{align}
(Note that this observation immediately reduces Proposition~\ref{prop:glob-injective} to Proposition~\ref{prop:loc-injective}.)

\smallskip\noindent
Condition (\ref{eqn:cond-inj}) says that   there exist $\left\lbrace \alpha_{ij} \mid \left[i,j \right) \in \mathcal{T}_1 \right\rbrace$,$\left\lbrace \beta_{ij} \mid \left[i,j \right) \in \mathcal{T}_2 \right\rbrace$ such that,
\begin{align}\label{eq:fundamental_lemma}
	\left( \forall \left[i,j\right) \in \mathcal{T}_1 \right) \alpha_{ij}&>0 \\
	\left( \forall \left[i,j\right) \in \mathcal{T}_2 \right) \beta_{ij}&>0 \\
    \sum_{ \left[i,j\right) \in \mathcal{T}_1} \alpha_{ij} &=1=\sum_{ \left[i,j\right) \in \mathcal{T}_2} \beta_{ij} \\
    \sum_{ \left[i,j\right) \in \mathcal{T}_1} \alpha_{ij}\left( e_i-e_j\right) &= \sum_{ \left[i,j\right) \in \mathcal{T}_2} \beta_{ij} \left( e_i-e_j \right) \label{eqn:spec-1}
\end{align}
Our objective is to show that conditions (\ref{eq:fundamental_lemma})--(\ref{eqn:spec-1}) imply $\mathcal{T}_1=\mathcal{T}_2$  and $\alpha_{i,j} = \beta_{i,j}$ for each interval $[i, j)\in \mathcal{T}_1=\mathcal{T}_2$.

 We begin with the observation that Proposition~\ref{prop:glob-injective} is trivially true for $n=3$. (In this case both $\partial((W_3)^\circ)$ and $\partial(Root_3)$ are boundaries of a hexagon.) This is sufficient to start an inductive proof. However note that we already know (Figure~\ref{fig-Root4-1}, Example~\ref{exam:3}, and Remark~\ref{rem:korisno}) that Proposition~\ref{prop:glob-injective} is also true in the case $n=4$.

\medskip
The proof is continued by induction on the parameter  $\nu :=\vert \mathcal{T}_1\vert + \vert \mathcal{T}_2\vert + n$. More precisely, we show that if there is a counterexample $(\mathcal{T}_1, \mathcal{T}_2)$ on $[n]$ then there is a counterexample $(\mathcal{T}_1', \mathcal{T}_2')$ on $[n']$ such that $\nu' = \vert \mathcal{T}_1'\vert + \vert \mathcal{T}_2'\vert + n' < \vert \mathcal{T}_1\vert + \vert \mathcal{T}_2\vert + n = \nu$.

\medskip\noindent
{Step~1:}  Without loss of generality we are allowed to assume that,
 \begin{equation} \label{eqn:so-si}
 s(\mathcal{T}_1)=s(\mathcal{T}_2)=I \mbox{ {\rm and} }  t(\mathcal{T}_1)=t(\mathcal{T}_2)=J.
 \end{equation}
 Indeed,   it follows from equation  (\ref{eqn:spec-1}) that $I$ (respectively $J$) collects the indices $i$ (respectively the indices $j$) where $e_i$ appears with a positive coefficient ($e_j$ appears with a negative coefficient).

 Moreover, we assume that,
 \begin{equation}\label{eqn:use-all}
   I \cup J = [n].
 \end{equation}
Otherwise, there exists an element $i\in [n]$ which is neither source nor terminal point of an interval in  $\mathcal{T}_1\cup \mathcal{T}_2$. In this case the vertex $i$ can be deleted and $n$ can be replaced by a smaller number $n'$.

 \medskip\noindent
 {Step~2:}  Let us assume that either $I$ or $J$ contains two consecutive elements, for example let $\{i, i+1\}\subset I$ for some $i\in [n]$. The proof in the case $\{i, i+1\}\subset J$ is similar (alternatively we can apply the automorphism $\beta$ from Proposition~\ref{prop:auto} which reverses the orientation of arcs).

 By cyclic relabelling, in other words by applying repeatedly the automorphism $\alpha$ from Proposition~\ref{prop:auto}, we may assume that $i=n-1$ and $i+1=n$.

 \bigskip
 Let $L_n  :  \mathbb{R}^n\rightarrow \mathbb{R}^{n-1}$ be the linear map such that $L_n(e_j) = e_j$ for $j=1,\dots, n-1$ and $L_n(e_n) = e_{n-1}$. On applying the map $L_n$ to both sides of the equality (\ref{eqn:spec-1}) we obtain a new relation,

 \begin{equation}\label{eqn:spec-2}
   \sum_{ \left[i,j\right) \in \mathcal{T}_1'} \alpha_{ij}'\left( e_i-e_j\right) = \sum_{ \left[i,j\right) \in \mathcal{T}_2'} \beta_{ij}' \left( e_i-e_j \right).
 \end{equation}
 For better combinatorial understanding of the relation (\ref{eqn:spec-2}), we note that a combinatorial-geometric counterpart of the map $L_n$ is the operation of collapsing  the interval $\left[ n-1,n \right]$ (in the circle $S^1$) to the point $n-1$.

\begin{figure}[htb]
	\centering
    \begin{tikzpicture}[scale=1]
    	\def \N {12} 
        \def \R {3} 
        \def \E {10/0/1,10/1/1,10/5/0,9/1/0,9/5/1,9/6/1,9/7/1}

        \pgfmathsetmacro{\ang}{360/\N}	
        \pgfmathsetmacro{\tmp}{\N-1};

		\draw (0,0) circle(\R);

        \foreach \x in {0,...,\tmp} {
           	\node[draw,circle,fill=white,inner sep=0pt, minimum size=6pt] (n\x) at ({\x*\ang}:{\R}) {};
      	}
        \node () at ({10*\ang}:{\R+.5}) {$n$};
        \node () at ({9*\ang}:{\R+.5}) {$n-1$};
        \node () at ({6*\ang}:{\R+.5}) {$B$};
        \node () at ({0.5*\ang}:{\R+.5}) {$A$};

       \foreach \i/\j/\t in \E {
        	\node[draw,circle,fill=black,inner sep=0pt, minimum size=6pt] (n\i) at ({\i*\ang}:{\R}) {};
            \pgfmathparse{abs(\j - \i) < (\N/2) ? int(1) : int(0)};
    		\ifnum\pgfmathresult=1
            	\pgfmathparse{\t > 0.5 ? int(1) : int(0)};
                \ifnum\pgfmathresult=1
					\draw[very thick,->-] (n\i) to[bend right] (n\j);
                \else
					\draw[very thick,dashed,->-] (n\i) to[bend right] (n\j);
                \fi
	   		\else
            	\pgfmathparse{\t > 0.5 ? int(1) : int(0)};
                \ifnum\pgfmathresult=1
					\draw[very thick,->-] (n\i) to[bend left] (n\j);
                \else
					\draw[very thick,dashed,->-] (n\i) to[bend left] (n\j);
                \fi
    		\fi
        	
        }
	\end{tikzpicture}
    \caption{Admissible family of arcs before the collapse of the interval $[n-1,n]$.}
	\label{fig-circ-1}
\end{figure}
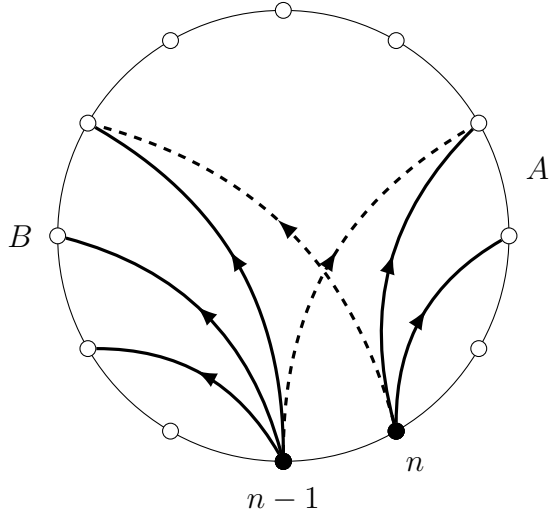

\smallskip
It is not difficult to describe the effect of the collapsing operation  ($CO$) on the admissible families $\mathcal{T}_1$ and $\mathcal{T}_2$ satisfying the condition $\{n-1, n\}\subset s(\mathcal{T}_1) = s(\mathcal{T}_2)$.

\begin{lema}\label{lema:co}
The collapsing operation $CO$  sends each admissible family $\mathcal{T}$ with the property $\{n-1, n\}\subset s(\mathcal{T})$ to an admissible family $\mathcal{T}'$ on the vertex set $[n-1] = \{1,\dots, n-1\}$. Moreover, under this condition, $CO([i,j)) = [i,j)$ if $i\neq n$, while $CO([n,j)) = [n-1,j)$.
\end{lema}
Each admissible family $\mathcal{T}$, satisfying the condition $\{n-1, n\}\subset s(\mathcal{T})$,
has a decomposition $\mathcal{T} = \mathcal{T}^a\uplus \mathcal{T}^b\uplus \mathcal{T}^c$ where,
 \begin{equation}\label{eqn:T-decom}
    \mathcal{T}^a = \{[i, j)\in \mathcal{T} \mid i = n\},\,   \mathcal{T}^b = \{[i, j)\in \mathcal{T} \mid i = n-1\}, \,
    \mathcal{T}^c = \mathcal{T}\setminus (\mathcal{T}^a\cup \mathcal{T}^b).
 \end{equation}
Let $A=A(\mathcal{T}) := \{j\in [n] \mid [n, j)\in \mathcal{T}\}$ and $B=B(\mathcal{T}) := \{j\in [n] \mid [n-1, j)\in \mathcal{T}\}$. Note that the sets $A(\mathcal{T})$ and $B(\mathcal{T})$ are either disjoint or have exactly one point in common. (Figure~\ref{fig-circ-1} shows how the common point arises as the end-point of one of the dotted arcs.)

\smallskip
It follows from Lemma~\ref{lema:co} that the admissible family $\mathcal{T}' := CO(\mathcal{T})$ admits the decomposition,
\begin{equation}\label{eqn:co-decomp}
\mathcal{T}' = CO(\mathcal{T}) = \mathcal{T}^c\uplus \mathcal{T}^{ab}
\end{equation}
where $\mathcal{T}^{ab} := \{[n-1, j) \mid j\in A(\mathcal{T})\cup B(\mathcal{T})\}$. This analysis and a  comparison of equalities (\ref{eqn:spec-1}) and (\ref{eqn:spec-2}) lead to the following observations.

\begin{enumerate}
 \item  By the induction hypothesis, the equality (\ref{eqn:spec-2}) leads to the conclusion that $\mathcal{T}_1' = \mathcal{T}_2'$ and $\alpha_{i,j}' = \beta_{i,j}'$ for  each pair $(i,j)$ such that $[i,j)\in \mathcal{T}_1' = \mathcal{T}_2'$;
 \item $\mathcal{T}_1' = \mathcal{T}_2'$ together with (\ref{eqn:co-decomp}) implies  $T^c_1 = T^c_2$ and $\alpha'_{i,j} = \alpha_{i,j} = \beta_{i,j} = \beta_{i,j}'$ for each pair $(i,j)$ such that $[i,j)\in \mathcal{T}^c_1 = \mathcal{T}^c_2$. By canceling these terms in (\ref{eqn:spec-1}) we obtain the following equality,
     \begin{equation}\label{eqn:spec-3}
       \sum_{ \left[i,j\right) \in \mathcal{T}_1^a\cup T_1^b} \alpha_{ij}\left( e_i-e_j\right) = \sum_{ \left[i,j\right) \in \mathcal{T}_2^a\cup T_2^b} \beta_{ij} \left( e_i-e_j \right)
     \end{equation}
  \item  This cancelation in (\ref{eqn:spec-3}) can be continued. The only indices, end-points of oriented arcs, that remain unaffected by the cancelation belong to the set $W:=(A(\mathcal{T}_1)\cap B(\mathcal{T}_1))\cup (A(\mathcal{T}_2)\cap B(\mathcal{T}_2))$.
   \item The only case that requires further analysis is the case when $W =\{i, j\}$ has precisely two elements (Figure~\ref{fig-circ-1}). In this case we obtain a contradiction by the argument used in Example~\ref{exam:3} (Remark~\ref{rem:korisno} and Figure~\ref{fig-circ-2}).
 \end{enumerate}

\medskip\noindent
Step~3: In this step we handle the only remaining case where neither $s(\mathcal{T}_1) = s(\mathcal{T}_2)=I$ nor  $t(\mathcal{T}_1) = t(\mathcal{T}_2)=J$ have consecutive elements. In this case there must exist two consecutive indices $\{i-1, i\}$ such that $i-1\in I$ and $i\in J$. Again, by the cyclic re-enumeration, we can assume that $i=n$ (Figure~\ref{fig-circ-3-4}, cases (a) and (b)).

\begin{figure}[h]
	\centering
    \begin{subfigure}[b]{0.4\textwidth}
    \begin{tikzpicture}[scale=0.9]
    	\def \N {12} 
        \def \R {3} 
        \def \E {9/5/1,9/6/1,9/7/1}

        \pgfmathsetmacro{\ang}{360/\N}	
        \pgfmathsetmacro{\tmp}{\N-1};

		\draw (0,0) circle(\R);

        \foreach \x in {0,...,\tmp} {
           	\node[draw,circle,fill=white,inner sep=0pt, minimum size=6pt] (n\x) at ({\x*\ang}:{\R}) {};
      	}
        \node () at ({10*\ang}:{\R+.5}) {$n$};
        \node () at ({9*\ang}:{\R+.5}) {$n-1$};
        \node () at ({6*\ang}:{\R+.5}) {$A$};

       \foreach \i/\j/\t in \E {
        	\node[draw,circle,fill=black,inner sep=0pt, minimum size=6pt] (n\i) at ({\i*\ang}:{\R}) {};
            \pgfmathparse{abs(\j - \i) < (\N/2) ? int(1) : int(0)};
    		\ifnum\pgfmathresult=1
            	\pgfmathparse{\t > 0.5 ? int(1) : int(0)};
                \ifnum\pgfmathresult=1
					\draw[very thick,->-] (n\i) to[bend right] (n\j);
                \else
					\draw[very thick,dashed,->-] (n\i) to[bend right] (n\j);
                \fi
	   		\else
            	\pgfmathparse{\t > 0.5 ? int(1) : int(0)};
                \ifnum\pgfmathresult=1
					\draw[very thick,->-] (n\i) to[bend left] (n\j);
                \else
					\draw[very thick,dashed,->-] (n\i) to[bend left] (n\j);
                \fi
    		\fi
        	
                \draw[very thick,dashed,->-] (n9) to[bend left] (n10);
                }
	\end{tikzpicture}
    \caption{}
    \end{subfigure}
    \begin{subfigure}[b]{0.1\textwidth}
    \end{subfigure}
    \begin{subfigure}[b]{0.4\textwidth}
    \begin{tikzpicture}[scale=0.9]
    	\def \N {12} 
        \def \R {3} 
        \def \E {10/0/1,10/1/1}

        \pgfmathsetmacro{\ang}{360/\N}	
        \pgfmathsetmacro{\tmp}{\N-1};

		\draw (0,0) circle(\R);

        \foreach \x in {0,...,\tmp} {
           	\node[draw,circle,fill=white,inner sep=0pt, minimum size=6pt] (n\x) at ({\x*\ang}:{\R}) {};
      	}
        \node () at ({10*\ang}:{\R+.5}) {$n$};
        \node () at ({9*\ang}:{\R+.5}) {$n-1$};
        \node () at ({0.5*\ang}:{\R+.5}) {$B$};

    \draw[very thick,dashed,->-] (n9) to[bend left] (n10);

	\draw[very thick,->-] (n0) to[bend right] (n10);
    \draw[very thick,->-] (n1) to[bend right] (n10);
    
\node[draw,circle,fill=black,inner sep=0pt, minimum size=6pt] (n0) at ({0*\ang}:{\R}) {};
\node[draw,circle,fill=black,inner sep=0pt, minimum size=6pt] (n1) at ({1*\ang}:{\R}) {};
\node[draw,circle,fill=black,inner sep=0pt, minimum size=6pt] (n9) at ({9*\ang}:{\R}) {};

	\end{tikzpicture}
    \caption{}
    \end{subfigure}
    \caption{}
	\label{fig-circ-3-4}
\end{figure}

For an admissible family $\mathcal{T}$, satisfying these conditions, there is a decomposition,
$\mathcal{T} = \mathcal{T}^a\cup \mathcal{T}^b\cup \mathcal{T}^c$, similar to (\ref{eqn:T-decom}), where $\mathcal{T}^a = \{[i, j)\in \mathcal{T} \mid j = n\},   \mathcal{T}^b = \{[i, j)\in \mathcal{T} \mid i = n-1\}$, and $\mathcal{T}^c = \mathcal{T}\setminus (\mathcal{T}^a\cup \mathcal{T}^b)$. Note that $\mathcal{T}^a\cap \mathcal{T}^b$ is either empty or $\mathcal{T}^a\cap \mathcal{T}^b =\{[n-1,n)\}$. It follows that $\mathcal{T}\setminus\{[n-1,n)\} = \hat{\mathcal{T}}^a\uplus \hat{\mathcal{T}}^b\uplus \mathcal{T}^c$, where $\hat{\mathcal{T}}^a := \mathcal{T}^a\setminus\{[n-1,n)\}$ and $\hat{\mathcal{T}}^b := \mathcal{T}^b\setminus\{[n-1,n)\}$.

\medskip
The key observation is that if $J_1 = [n-1, j_1)\in \hat{\mathcal{T}}^a$ and $I_1 = [i_1, n)\in \hat{\mathcal{T}}^b$, then intervals $I_1$ and $J_1$ intersect but cannot be compatible in the sense of Definition~\ref{def:admissible}. It immediately follows that either $\hat{T}^a = \emptyset$ or $\hat{\mathcal{T}}^b = \emptyset$.

\begin{lema}\label{lema:AB}
  If $\mathcal{T}_1$ and $\mathcal{T}_2$ are (Step~3) admissible families satisfying conditions (\ref{eq:fundamental_lemma})-(\ref{eqn:spec-1}) then either,
  \begin{enumerate}
    \item[{\rm (A)}]\qquad  $\hat{\mathcal{T}}_1^a = \hat{\mathcal{T}}_2^a = \emptyset$\quad  and \quad $\mathcal{T}_i^b = \hat{\mathcal{T}}_i^b \cup \{[n-1,n)\}$\,  for $i=1,2$ , or
    \item[{\rm (B)}]  \qquad $\hat{\mathcal{T}}_1^b = \hat{\mathcal{T}}_2^b = \emptyset$ \quad and \quad $\mathcal{T}_i^a = \hat{\mathcal{T}}_i^a \cup \{[n-1,n)\}$\,  for $i=1,2$.
  \end{enumerate}
  \end{lema}
For illustration, the case $\hat{\mathcal{T}}_1^a = \hat{\mathcal{T}}_2^b = \emptyset$ is ruled out by the following reasoning. The only way to satisfy the condition (\ref{eqn:spec-1}) is to have $[n-1,n)\in \mathcal{T}_1^b\cap \mathcal{T}_2^a$  (and $\alpha_{n-1,n} \neq 0 \neq \beta_{n-1,n}$). This is not possible, however, since by comparing the coefficients near $e_{n-1}$ and $e_n$ we obtain the following contradictory equalities,
\begin{equation}\label{eqn:contradictory}
  \alpha_{n-1,n} + \sum_{[n-1,j)\in \mathcal{T}_1^b} \alpha_{n-1,j}= \beta_{n-1,n}  \qquad \alpha_{n-1,n} = \beta_{n-1,n} + \sum_{[i,n)\in \mathcal{T}_2^a} \beta_{i,n}.
\end{equation}

The proof (Step~3) is continued by observing that in both cases of Lemma~\ref{lema:AB} the coefficients $\alpha_{n-1,n}$ and $\beta_{n-1,n}$ must be equal, and the corresponding terms in (\ref{eqn:spec-1}) can be cancelled out.
The proof (Step~3) is finished by applying the collapsing operator $CO$, which collapses the interval $[n-1,n]$ either to the left end-point $(n-1)$ (corresponding to the case A of Lemma~\ref{lema:AB}) or to the right end-point $n$  (corresponding to the case B of Lemma~\ref{lema:AB}). The analysis is similar to the collapsing procedure described in Step~2 so we omit the details.

\medskip
This completes the proof of Proposition~\ref{prop:glob-injective}. \hfill $\square$

\subsection{Surjectivity of the map $\phi_n$ }
\label{sec:surjectivity}

  We already know (Proposition~\ref{prop:glob-injective}) that $\phi_n$ is injective. Let us show that it is surjective as well.

\smallskip By Proposition~\ref{prop:glob-injective} the map $\phi_n$ induces an isomorphism in homology, i.e. the degree ${\rm deg}(\phi_n)$ is either $+1$ or $-1$. Therefore it must be an epimorphism since otherwise it would be homotopic to a constant map.   \hfill $\square$

\section{Cyclohedron and the root polytope II}
\label{sec:revisited}

 In this section we give a geometric interpretation of the map $\phi_n : \partial(W_n^\circ)\rightarrow \partial(Root_n)$, introduced in Definition~\ref{def:canonical-map}. The key observation is that the dual $(Root_n)^\circ$ of the root polytope belongs to the {\em  irredundant part} of the   {\em face deformation cone} \cite[Section~15]{prw} of the cyclohedron $W_n$. A more precise statement says that the pair $(W_n, (Root_n)^\circ)$  may be interpreted as a couple of polytopes $(\Delta_{\mathcal{\widehat{F}}}, \Delta_{\mathcal{F}})$, where $\mathcal{\widehat{F}}$ is a building set and $\mathcal{F}$ its irredundant basis in the sense of Definition~\ref{def:irredundant}. For an introduction into the theory of {\em nestohedra, building sets, etc.}, the reader is referred to \cite{bp,po,fs}.

\subsection{Building set of the cyclohedron $W_n$}
\label{sec:building}

It is well-known, see \cite[Section~3]{fs}, \cite{po}, \cite{devadoss}, or \cite[Section~1.5.]{bp}, that the cyclohedron $W_n$ is a nestohedron (graph associahedron), so it has a Minkowski sum decomposition,
\begin{equation}\label{eqn:cycl-Mink-1}
  W_n = \Delta_{\mathcal{\widehat{F}}} :=\sum_{F\in \mathcal{\widehat{F}}} \Delta_F,
\end{equation}
where $\mathcal{\widehat{F}}\subset 2^{[n]}\setminus\{\emptyset\}$ is the associated {\em building set} \cite{fs, po} and,
$$
\Delta_F\subset \Delta_{[n]}:= {\rm Conv}\{e_i\}_{i=1}^n\subset \mathbb{R}^n,
$$
is the simplex spanned by $F\subset [n]$. The family $\mathcal{\widehat{F}}$ is in the case of $W_n$ identified as the collection $\mathcal{\widehat{F}} = Con(C_n)$ of all connected subsets in the cycle graph $C_n$ with $n$-vertices. Note that a set  $Z\subset [n]$ is connected if $Z$ is either a cyclic interval or $Z=[n]$.

\bigskip
The Minkowski sum $\Delta_{\mathcal{F}} :=\sum_{F\in \mathcal{F}} \Delta_F$ is defined for any family (hypergraph) $\mathcal{F}\subset 2^{[n]}\setminus\{\emptyset\}$, however it is not necessarily a simple polytope, unless $\mathcal{F}$ is a building set. For this reason it is interesting to compare $\Delta_{\mathcal{F}}$ and $\Delta_{\mathcal{\widehat{F}}}$ where $\mathcal{\widehat{F}}$ is the {\em building closure} of $\mathcal{F}$.

\begin{defin}\label{def:building-closure}
  A family $\mathcal{\widehat{F}}\supset \mathcal{F}$ is the {\em building closure} of a hypergraph $\mathcal{F}$ if $\mathcal{\widehat{F}}$ is the unique minimal building set which contains $\mathcal{F}$. In this case we also say that $\mathcal{{F}}$ is a {\em building basis} of the building set $\mathcal{\widehat{F}}$.
\end{defin}

\begin{defin}\label{def:build-restr-del}
If $\mathcal{F}\subset 2^{[n]}\setminus\{\emptyset\}$ is a hypergraph and $X\subset [n]$ then,
\[
\mathcal{F}_X := \{F\in \mathcal{F} \mid F\subset X\} \quad {\rm and}  \quad \mathcal{F}^X := \{F\setminus X \mid F\in \mathcal{F},  F\setminus X \neq\emptyset\}.
\]
The family $\mathcal{F}_X$ is referred to as the {\em restriction} of $\mathcal{F}$ to $X$, while $\mathcal{F}^X$ is obtained from $\mathcal{F}$ by {\em deletion} of the set $X$.
\end{defin}

 For each $X\subset [n]$ let $\phi_X : \mathbb{R}^n \rightarrow \mathbb{R}$ be the linear form $\phi_X(x) = \sum_{i\in X}~x_i$ (for example $\phi_{[n]}(x) = x_1+x_2\dots+x_n$).  The cardinality of a family $\mathcal{F}$ is denoted by $\vert \mathcal{F}\vert$.

 \medskip
 The following proposition, see \cite[Proposition~7.5.]{po} or \cite[Proposition~3.12.]{fs}, provides a description of $\Delta_{\mathcal{F}}$ in terms of linear (in)equalities.

 \begin{prop}\label{prop:3.12-fs}
 Suppose that $\mathcal{\widehat{F}}$ is the building closure of a family
 $\mathcal{F}\subset 2^{[n]}\setminus\{\emptyset\}$. Then,
 \begin{equation}\label{eqn:irred}
 \Delta_{\mathcal{F}} =  \left\{ x\in \mathbb{R}^n \mid \phi_{[n]}(x) = \vert \mathcal{F}\vert \,  \mbox{ {\rm and for each} }\, X\in \mathcal{\widehat{F}},   \,  \phi_X(x) \geqslant \vert \mathcal{F}_X\vert \right\}.
 \end{equation}
 Moreover, the face of $\Delta_{\mathcal{F}}$ where $\phi_X$ attains its minimum is isomorphic to the Minkowski sum,
 \begin{equation}\label{eqn:attain-min}
 {\rm Face}_{\phi_X}(\Delta_{\mathcal{F}})  \cong \Delta_{\mathcal{F}_X} +   \Delta_{\mathcal{F}^X} \, .
 \end{equation}
   \end{prop}

\begin{defin}\label{def:irredundant}
  We say that a hypergraph $\mathcal{F}\subset 2^{[n]}\setminus\{\emptyset\}$ is {\em tight} if all inequalities in (\ref{eqn:irred})  are essential (irredundant), where $\mathcal{\widehat{F}}$ is the building closure of $\mathcal{F}$. We also say that $\mathcal{F}$ is a tight or {\em irredundant basis} of the building set $\mathcal{\widehat{F}}$.
\end{defin}

 The following criterion for tightness of $\mathcal{F}$ is easily deduced from Proposition~\ref{prop:3.12-fs}.

 \begin{prop}\label{prop:irredund}
   Let $\mathcal{F}\subset 2^{[n]}\setminus\{\emptyset\}$ be a hypergraph and let $\mathcal{\widehat{F}}$ be its building closure. Then $\mathcal{F}$ is tight if and only for each $X\in \mathcal{\widehat{F}}$,
   \begin{enumerate}
     \item[{\rm (1)}] $\mathcal{F}_X$ is connected as a hypergraph on $X$ and,
     \item[{\rm (2)}] $\mathcal{F}^X$ is connected as a hypergraph on $[n]\setminus X$.
   \end{enumerate}
 Actually, the first condition is automatically satisfied, as a consequence of the fact that  $\mathcal{\widehat{F}}$ is the building closure of $\mathcal{F}$.
 \end{prop}

\medskip\noindent
{\bf Proof:} By Proposition~\ref{prop:3.12-fs} (relation (\ref{eqn:attain-min}))
$\mathcal{F}$ is tight if and only if for each $X\in \mathcal{\widehat{F}}$,
 \[
     {\rm dim}(\Delta_{\mathcal{F}_X} + \Delta_{\mathcal{F}^X}) = {\rm dim}(\Delta_{\mathcal{F}}) -1 = n-2.
 \]

\smallskip
It is well known, see \cite[Proposition~1.5.2.]{bp} or \cite[Remark~3.11.]{fs}, that for each hypergraph $\mathcal{H}\subset 2^S\setminus\{\emptyset\}$ the dimension of the polytope $\Delta_{\mathcal{H}} = \sum_{H\in \mathcal{H}}~\Delta_H$ is $\vert S \vert -c$ where $c$ is the number of components of the hypergraph $\mathcal{H}$. (Recall that $x, y\in S$ are in the same connected component in the hypergraph $\mathcal{H}$ if there is a sequence of elements $x = z_1, z_2,\dots, z_k = y$ in $S$ such that each $\{z_i, z_{i+1}\}$  is contained in some `edge' of the hypergraph $\mathcal{H}$.)

\smallskip
Since $\mathcal{\widehat{F}}$ is the building closure of $\mathcal{F}$, ${\rm dim}(\Delta_{\mathcal{F}_X})= \vert X \vert -1$. Indeed, by the proof of \cite[Lemma~3.10]{fs} $X\in \mathcal{\widehat{F}}$ if and only if $X$ is a singleton or $\mathcal{F}_X$ is a connected hypergraph on $X$.

\smallskip
It follows that,
\[
{\rm dim}(\Delta_{\mathcal{F}_X} + \Delta_{\mathcal{F}^X}) =
{\rm dim}(\Delta_{\mathcal{F}_X}) + {\rm dim}(\Delta_{\mathcal{F}^X}) = n-2 = (\vert X\vert -1) + (n-\vert X\vert -1)
\]
if and only if ${\rm dim}(\Delta_{\mathcal{F}^X}) = n-\vert X\vert -1$. This equality is equivalent to the condition (2) in Proposition~\ref{prop:irredund}. \hfill $\square$

\medskip
As a consequence of Proposition~\ref{prop:irredund} we obtain the following result.

\begin{prop}\label{prop:ciklo-irred}
Let $C_n$ be the cycle graph on $n$ vertices and let {\rm (\ref{eqn:cycl-Mink-1})} be the associated graph associahedron representation of the cyclohedron $W_n$, where $\mathcal{\widehat{F}} = Con(C_n)$ is the building set of all $C_n$-connected subsets in $[n]$. Then,
\begin{equation}\label{eqn:1-2-3}
\mathcal{F} := \{\{1,2\}, \{2,3\}, \dots, \{n-1,n\}, \{n,1\}\}
\end{equation}
is a tight hypergraph on $[n]$, which is a tight (irredundant) basis of $\mathcal{\widehat{F}}$ in the sense of Definition~\ref{def:irredundant}. As a consequence all inequalities {\rm (\ref{eqn:irred})}, in the corresponding description of the polytope,
\begin{equation}\label{eqn:rombic-gen-irred}
  \Delta_{\mathcal{F}} = \sum_{i=1}^n~[e_i, e_{i+1}]
\end{equation}
are essential (irredundant).
\end{prop}

\medskip\noindent
{\bf Proof:} If $X=[i,j]^0 = \{i, i+1,\dots, j\}\subset [n] = {\rm Vert}(C_n)$ is an element in $\mathcal{\widehat{F}}= Con(C_n)$ then,
\[
       \mathcal{F}^X = \{ \{i-1\}, \{j+1\}\}\cup  \{[j+1,j+2]^0, [j+2, j+3]^0, \dots, [i-2, i-1]^0\}.
\]
This hypergraph is clearly connected on its set of vertices which verifies the condition (2) in Proposition~\ref{prop:irredund} and completes the proof of Proposition~\ref{prop:ciklo-irred}. \hfill $\square$

\subsection{Polar dual of the root polytope}\label{sec:dual}

\begin{defin}\label{def:delta-zono}
Let ${\Delta = \Delta_A = \rm Conv}(A) = {\rm Conv}\{a_0, a_1,\dots, a_m\}$ be a non-degenerate simplex with vertices in $A\subset \mathbb{R}^n$. The associated $\Delta$-zonotope $Zono(\Delta) = Zono(\Delta_A)$ is the Minkowski sum,
  \begin{equation}\label{eqn:delta-zono}
    Zono(\Delta_A) = [a, a_0]+[a, a_1]+\dots+[a, a_m],
  \end{equation}
 where $a := \frac{1}{m+1}(a_0+\dots+ a_m)$. If $\Delta_n = \Delta_{[n]} = {\rm Conv}(\{e_1,\dots, e_n\})$ is the simplex spanned by the orthonormal basis in $\mathbb{R}^n$, then the associated $\Delta$-zonotope $Zono(\Delta_n) = Zono(\Delta_{[n]})$ is referred to as the standard, $(n-1)$-dimensional $\Delta$-zonotope.
\end{defin}

\begin{figure}[hbt]
\centering 
\includegraphics[scale=0.45]{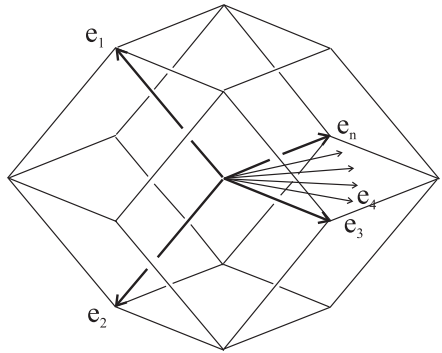}
 \caption{$\Delta$-zonotope $Zono(\Delta_n)$ as a generalized rhombic dodecahedron.}
 \label{fig:R12}
\end{figure}
\begin{exam}

  Figure~\ref{fig:R12} depicts the standard $(n-1)$-dimensional $\Delta$-zonotope.
  In the special case $n=4$ one obtains the rhombic dodecahedron.
\end{exam}

\begin{defin}\label{def:gen-root}
 The generalized root polytope associated to a simplex $\Delta_A = {\rm Conv}(A) = {\rm Conv}\{a_1,\dots, a_n\}$ is the polytope,
 \begin{equation}\label{ref:gen-root-1}
     Root(\Delta_A) = {\rm Conv}\{a_i-a_j \mid 1\leqslant i\neq j\leqslant n\}.
 \end{equation}
 If  $W$ is an affine map such that $a_i = W(e_i)$ for each $i\in [n]$ then,
 \begin{equation}\label{eqn:gen-root-2}
   Root(W(\Delta_n)) = W(Root(\Delta_n)) = W(Root_n).
 \end{equation}
 \end{defin}

\medskip
The root polytope $Root_n$ is a subset of the hyperplane $H_0 =\{x\in \mathbb{R}^n \mid x_1+\dots+ x_n =0 \}$ while $Zono(\Delta_n) \subset H_n =\{x\in \mathbb{R}^n \mid x_1+\dots+ x_n =n \}$. In the following proposition we claim that the $\Delta$-zonotope $Zono(\Delta_n^0)$, obtained by translating  $Zono(\Delta_n)$ to $H_0$, is precisely the polar dual of the polytope $Root_n$.

\begin{prop}\label{prop:dual}
  The root polytope $Root_n = {\rm Conv}\{e_i-e_j \mid 1\leqslant i\neq j\leqslant n\}\subset H_0$ is the dual (inside $H_0$) of the $\Delta$-zonotope $Zono(\Delta_{n}^0)$ where $\Delta_{n}^0 = \Delta_{n} - \frac{1}{n}(e_1+\dots+e_n)$,
\begin{equation}\label{eqn:root-polar}
  (Root_n)^\circ  =  Zono(\Delta_{n}^0).
\end{equation}
\end{prop}

\medskip\noindent
{\bf Proof:} The proof is an elementary exercise in the concept of duality (see \cite[Proposition~7]{ziv15}).  Let $\hat{e}_i = e_i - \frac{1}{n}(e_1+\dots + e_n)\in H_0$. It is sufficient to observe that the dual of the root polytope is,
\begin{equation}\label{eqn:root-dual}
  (Root_n)^\circ = \{x\in H_0 \mid   \vert x_i-x_j\vert \leqslant 1,\,  \mbox{ {\rm for} } 1\leqslant i \leqslant n)\},
\end{equation}
while the two supporting hyperplanes of $Zono(\Delta_n^0)$, parallel to $\mathcal{L}_{i,j} = {\rm span}\{\hat{e}_k \mid i\neq k\neq j\}$ (Figure~\ref{fig:R12}) have equations,
\[
     x_i - x_j =1  \quad   \mbox{\rm {and} } \quad  x_j-x_i =1.
\]
\hfill $\square$

\begin{lema}\label{lema:affine}
  Suppose that $K^\circ$ is the polar dual of a convex body $K$.  If $A : \mathbb{R}^d\rightarrow \mathbb{R}^d$ is a non-singular linear map then,
  \begin{equation}\label{eqn:affine}
         (A(K))^\circ = B(K^\circ)
  \end{equation}
 where  $B = (A^\ast)^{-1}$. In particular if $A = (A^\ast)^{-1}$ is an orthogonal transformation and $\mu\neq 0$ then, $(A(K))^\circ = A(K^\circ)$ and $(\mu K)^\circ = (1/\mu) K^\circ$.
\end{lema}

\medskip\noindent
{\bf Proof:}  If $y\in  (A(K))^\circ$ then by definition,
\begin{equation}\label{eqn:polarity-1}
    \langle y, Ax \rangle  = \langle A^\ast y, x \rangle \leqslant 1 \quad \mbox{\rm {for each}} \, x\in K.
\end{equation}

\medskip
The following extension of Proposition~\ref{prop:dual} is recorded for the future reference.

\begin{prop}\label{prop:dual-gen}
 Let $H_0\subset \mathbb{R}^n$ be the subspace spanned by $\hat{e}_i = e_i - \frac{1}{n}(e_1+\dots + e_n)$ and let $\Delta_n^0 = {\rm Conv}\{\hat{e}_1, \dots, \hat{e}_n\}$.  Let $A  : H_0 \rightarrow H_0$ be a non-singular linear map and let $B := (A^\ast)^{-1}$. Then,
 \begin{equation}\label{eqn:dual-gen}
    (Root(A(\Delta_n^0)))^\circ =  Zono(B(\Delta_n^0))
 \end{equation}
 \end{prop}

\smallskip\noindent
{\bf Proof:} By definition  $Root(A(\Delta_n^0)) = A(Root(\Delta_n^0))$ and $Zono(B(\Delta_n^0)) = B(Zono(\Delta_n^0))$. It follows from Proposition~\ref{prop:dual} and Lemma~\ref{lema:affine} that,
\[
   (Root(A(\Delta_n^0)))^\circ = (A(Root(\Delta_n^0)))^\circ = B((Root(\Delta_n^0))^\circ) = B(Zono(\Delta_n^0)) =  Zono(B(\Delta_n^0)).
\]

\subsection{$W_n^\circ$ as a Kantorovich-Rubinstein polytope}

There is a canonical isomorphism of vector spaces $H_\nu = \{x\in \mathbb{R}^n \mid x_1+\dots + x_n =\nu\}$ and the quotient space  $\mathbb{R}^n/\mathbb{R}e$, where $e = \frac{1}{n}(e_1+\dots + e_n)$, which induces a canonical isomorphism between $H_0$ and $H_\nu = \{x\in \mathbb{R}^n \mid x_1+\dots + x_n =\nu\}$ for each $\nu\in \mathbb{R}$.

\smallskip
The canonical isomorphism between $H_1$ and $H_0$ sends $e_i$ to $\hat{e}_i = e_i - \frac{1}{n}(e_1+\dots + e_n)$ and $\Delta_n = {\rm Conv}\{{e}_1, \dots, {e}_n\}$ to $\Delta_n^0 = {\rm Conv}\{\hat{e}_1, \dots, \hat{e}_n\}$.

\smallskip
The canonical isomorphism between $H_n$ and $H_0$ identifies the polytope $\Delta_{\mathcal{F}}$, introduced in  Section~\ref{sec:building} (equation  (\ref{eqn:rombic-gen-irred})), with the polytope,
\begin{equation}\label{eqn:link}
  \Delta_{\mathcal{F}} = \sum_{i=1}^n~[\hat{e}_i, \hat{e}_{i+1}]  = \sum_{i=1}^n~(\hat{e}_i+[0,  \hat{e}_{i+1}-\hat{e}_i]) =  \sum_{i=1}^n~[0, \hat{e}_{i+1}-\hat{e}_i] = Zono(B(\Delta_n^0))
\end{equation}
where $e_{n+1}:= e_1$ and $B : H_0\rightarrow H_0$ is the linear map defined by $B(\hat{e}_i) = b_i = \hat{e}_{i+1} - \hat{e}_i$.
In other words the polytope $\Delta_{\mathcal{F}}$  (associated to the irredundant basis (\ref{eqn:1-2-3}) of the building set  $\mathcal{\widehat{F}} = Con(C_n)$) is a $\Delta$-zonotope (generalized rhombic dodecahedron) $Zono(B(\Delta_n^0))$.

\medskip
The dual of $Zono(B(\Delta_n^0))$  is by Proposition~\ref{prop:dual-gen}  a root polytope, \begin{equation}\label{eqn:nov-root}
  Root_n^a := {\rm Conv}\{a_i - a_j \mid 1\leqslant i\neq j\leqslant n\} = Root(A(\Delta_n^0))
\end{equation}
where the vectors $\{a_i\}_{i=1}^n\subset H_0$ are defined by $a_i = A(\hat{e_i})$ and  $A : H_0\rightarrow H_0$ is the linear map such that $B = (A^\ast)^{-1}$ ($A = (B^\ast)^{-1}$ ).

\medskip
Summarizing, we record for the future reference the following proposition,

\begin{prop}\label{prop:nova-baza}
There exist vectors $\{a_i\}_{i=1}^n\subset H_0$, such that $a_1+\dots+ a_n=0$ and ${\rm Span}(\{a_i\}_{i=1}^n) = H_0$, which have the property that the dual of the root polytope $Root_n^a$ (defined by (\ref{eqn:nov-root})) is the polytope $\Delta_{\mathcal{F}}$ (defined by (\ref{eqn:link})).
\end{prop}

\medskip The following theorem is the main result of Section~\ref{sec:revisited}.

\begin{theorem}\label{thm:glavna-B}
  There exists a quasi-metric (asymmetric distance function) $\rho$ on the set $[n]$  such that the associated Kantorovich-Rubinstein polytope,
  \begin{equation}\label{eqn:KR-poly-2}
  KR(\rho) = {\rm Conv}\left\{\frac{e_i-e_j}{\rho(i,j)} \mid 1\leqslant i\neq j\leqslant n\right\}
  \end{equation}
  is {\em affinely isomorphic} to the dual $W_n^\circ$ of the cyclohedron $W_n$. Moreover the distance function $\rho$ satisfies a strict triangle inequality in the sense that,
  \[
  \rho(i,k) < \rho(i,j) + \rho(j,k) \, \mbox{ {\rm if} } \, i\neq j\neq k \neq i .
  \]
\end{theorem}

\smallskip\noindent
{\bf Proof:} Let $\{a_i\}_{i=1}^n\subset H_0$ be the collection of vectors described in  Proposition~\ref{prop:nova-baza} and let $a_{i,j} = a_i-a_j$ (for $1\leqslant i\neq j\leqslant n$) be the corresponding roots. In light of Proposition~\ref{prop:nova-baza} the polytope $\Delta_{\mathcal{F}}$ has the following description,
\begin{equation}\label{eqn:opis-1}
  \Delta_{\mathcal{F}}  = \{x\in H_0 \mid \langle a_{i,j}, x \rangle \leqslant 1  \, \mbox{ {\rm for each pair} }\, i\neq j \}.
\end{equation}
 All inequalities in (\ref{eqn:opis-1}) are irredundant. Moreover, the analysis from Section~\ref{sec:building} guarantees that there exist positive real numbers  $\{\alpha_{i,j} \mid 1\leqslant i\neq j\leqslant n \}$ such that,
 \begin{equation}\label{eqn:opis-2}
  W_n = \Delta_{\mathcal{\widehat{F}}}  = \{x\in H_0 \mid \langle a_{i,j}, x \rangle \leqslant \alpha_{i,j}  \, \mbox{ {\rm for each pair} }\, i\neq j \}.
\end{equation}
From here it immediately follows that,
\[
    W_n^\circ = {\rm Conv}\left\{\frac{a_{i,j}}{\alpha_{i,j}} \mid 1\leqslant i\neq j\leqslant n\right\}.
\]
Let us show that $\rho(i,j) := \alpha_{i,j}$ is a strict quasi-metric on $[n]$. Assume that there exist three, pairwise distinct, indices $i.j.k\in [n]$ such that $\rho(i,k)\geqslant \rho(i,j) + \rho(j,k)$. Then,
\[
    \frac{a_{i,k}}{\rho(i,k)}\in [0, a_{i,k}/(\rho(i,j) + \rho(j,k))]  \quad \mbox{{\rm and}} \quad  \frac{a_{i,k}}{\alpha_{i,k}}\in [0, a_{i,k}/(\alpha_{i,j} + \alpha_{j,k}))].
\]
In light of the obvious equality,
\[
    \frac{a_{i,k}}{\alpha_{i,j}+\alpha_{j,k}} = \frac{\alpha_{i,j}}{\alpha_{i,j}+\alpha_{j,k}}\left(\frac{a_{i,j}}{\alpha_{i,j}}\right) + \frac{\alpha_{j,k}}{\alpha_{i,j}+\alpha_{j,k}}\left(\frac{a_{j,k}}{\alpha_{j,k}}\right)
\]
we observe that if both inequalities,
\[
\langle a_{i,j}, x \rangle \leqslant \alpha_{i,j}  \quad \mbox{{\rm and}} \quad  \langle a_{j,k}, x \rangle \leqslant \alpha_{j,k}
\]
are satisfied then $\langle a_{i,k}, x \rangle \leqslant \alpha_{i,k}$. This is however in contradiction with non-redundancy of the last inequality in the representation (\ref{eqn:opis-2}). \hfill $\square$

\subsection{An explicit quasi-metric associated to a cyclohedron}
\label{sec:explicit}

\begin{defin}\label{def:height}
  Let $\mathcal{\widehat{F}}$ be the building closure of a hypergraph $\mathcal{F}\subset 2^{[n]}\setminus\{\emptyset\}$. The associated {\em `height function'} $h_{\mathcal{F}}  :  \mathcal{\widehat{F}} \rightarrow \mathbb{R}$
  is defined by, $$h_{\mathcal{F}}(X) = \frac{\vert\mathcal{F}\vert}{n} - \frac{\vert \mathcal{F}_X\vert}{\vert X\vert},$$
 so in particular $h_{\mathcal{F}}([n])=0$ for each hypergraph $\mathcal{F}$.
\end{defin}

\medskip
The inequalities (\ref{eqn:irred}), describing $\Delta_{F}$ as a subset of
$H_{\vert \mathcal{F}\vert} = \{x\in \mathbb{R}^n \mid \phi_{[n]}(x) = \vert \mathcal{F}\vert \}$ can be, with the help of the height function, rewritten as follows,
\begin{equation}\label{eqn:irred-2}
 \Delta_{\mathcal{F}} =  \left\{ x\in H_{\vert\mathcal{F}\vert} \mid  \mbox{ {\rm For each} }\, X\in \mathcal{\widehat{F}}\setminus\{[n]\},   \,  \frac{\phi_{[n]}(x)}{n}-\frac{\phi_X(x)}{\vert X\vert} \leqslant h_{\mathcal{F}}(X) \right\}.
 \end{equation}

In particular, if $\mathcal{F} = \mathcal{\widehat{F}}$ we obtain the representation,
\begin{equation}\label{eqn:irred-3}
 \Delta_{\mathcal{\widehat{F}}} =  \left\{ x\in H_{\vert\mathcal{\widehat{F}}\vert} \mid  \mbox{ {\rm For each} }\, X\in \mathcal{\widehat{F}}\setminus\{[n]\},   \,  \frac{\phi_{[n]}(x)}{n}-\frac{\phi_X(x)}{\vert X\vert} \leqslant h_{\mathcal{\widehat{F}}}(X) \right\}.
 \end{equation}
Assuming that $h_\mathcal{F}(X)\neq 0$ for each $X\in \mathcal{\widehat{F}}\setminus\{[n]\}$, let $A_X\in H_0$ be the vector defined by,
\begin{equation}\label{eqn:A-X}
  A_X := \frac{1}{h_\mathcal{F}(X)}\left( \frac{e}{n} - \frac{e_X}{\vert X\vert} \right),
\end{equation}
where $e_X := \sum_{i\in X}~e_i$ and $e = e_{[n]} = e_1+\dots + e_n$. It follows  that (\ref{eqn:irred-2}) and (\ref{eqn:irred-3}) can be rewritten as,
\begin{equation}\label{eqn:A-X-irred-2}
 \Delta_{\mathcal{F}} =  \left\{ x\in H_0 \mid  \mbox{ {\rm For each} }\, X\in \mathcal{\widehat{F}}\setminus\{[n]\},   \,  \langle A_X, x \rangle \leqslant 1 \right\},
 \end{equation}
\begin{equation}\label{eqn:A-X-irred-3}
 \Delta_{\widehat{\mathcal{F}}} =  \left\{ x\in H_0 \mid  \mbox{ {\rm For each} }\, X\in \mathcal{\widehat{F}}\setminus\{[n]\},   \,  \langle A_X, x \rangle \leqslant \frac{h_\mathcal{\widehat{F}}(X)}{h_\mathcal{F}(X)} \right\}.
 \end{equation}

Now we specialize to the case $\mathcal{F} := \{\{1,2\}, \{2,3\}, \dots, \{n-1,n\}, \{n,1\}\}$,    so the associated building closure  $\mathcal{\widehat{F}} = Con(C_n)$ is, as in Proposition~\ref{prop:ciklo-irred}, the building set of all $C_n$-connected subsets (circular intervals) in the circle graph $C_n$ . The corresponding height functions are shown in the following lemma.

\begin{lema}\label{lema:height}
If $X = [n]$ then $h_{\mathcal{F}}(X) = h_{\mathcal{\widehat{F}}}(X) = 0$. If $X\neq [n]$ then,
\[
 h_{\mathcal{F}}(X) =  \frac{1}{\vert X\vert}\quad  \mbox{ {\rm and} } \quad
h_{\mathcal{\widehat{F}}}(X) = \frac{n^2-n+1}{n} - \frac{\vert X\vert +1}{2}.
\]
\end{lema}

Recall (Section~\ref{sec:face}) that for $i, j\in {\rm Vert}(C_n) = [n]$, the associated (discrete) circular interval is $[i,j]^0 := [i,j]\cap [n]$. Similarly  (for $i\neq j$)  $[i,j)^0 := [i,j)\cap [n]$, so  $[i,j)^0 = [i,j-1]^0$ if $i\neq j$ and $[i,i)^0 = \emptyset$ . Define the {\em ``clock quasi-metric''} on ${\rm Vert}(C_n)$  by,
\begin{equation}\label{eqn:quasi-circ}
   d(i,j) := \vert [i,j)^0\vert = \vert [i,j]^0\vert -1.
\end{equation}

\begin{prop}\label{prop:quasi-metric-cyclo}
 Let $\rho$ be the quasi-metric on $[n] = {\rm Vert}(C_n)$ defined by,
 \begin{equation}\label{eqn:explicit-d}
   \rho(i,j) =  d(i,j)\frac{n^2-n+1}{n}- \frac{d(i,j)(d(i,j) +1)}{2}
 \end{equation}
 where $d$ is the clock quasi-metric on $[n]$.  Than the associated Kantorovich-Rubinstein polytope $KR(\rho)$ is affinely isomorphic to a polytope dual to the standard cyclohedron.
\end{prop}

\medskip\noindent
{\bf Proof:} By definition $\rho(i,j) = h_\mathcal{\widehat{F}}([i,j)^0)/h_\mathcal{F}([i,j)^0)$. Recall that equations (\ref{eqn:A-X-irred-2}) and (\ref{eqn:A-X-irred-3}) are nothing but a more explicit form of equations (\ref{eqn:opis-1}) and (\ref{eqn:opis-2}). It immediately follows that $\rho(i,j) = \alpha_{i,j}$ which by Theorem~\ref{thm:glavna-B} implies that $\rho$ is indeed a quasi-metric on $[n]$ such that the associated Lipschitz polytope $Lip(\rho)$  is a cyclohedron.  \hfill $\square$

\begin{remark}{\rm
  By a similar argument we already know that vectors $\{A_X \mid X\in \mathcal{\widehat{F}}\setminus\{[n]\}\}$,  described by equation (\ref{eqn:A-X}),  form a type A root system if $\mathcal{F} := \{\{1,2\}, \{2,3\}, \dots,  \{n,1\}\}$. This can be seen directly as follows. For $2\leq i\leq n$,  let $X_i := [1,i)^0$ and $A_i:= A_{X_i}$. Then there is a disjoint union,
  \[
      \{A_X \mid X\in \mathcal{\widehat{F}}\setminus\{[n]\}\} =  \{A_i\}_{i=2}^n \sqcup \{-A_i\}_{i=2}^n \sqcup \{A_i - A_j \mid 2 \leq i\neq j\leq n\}.
  \]
  This observation and a comparison of  (\ref{eqn:link}) and (\ref{eqn:A-X-irred-2}) provide an alternative proof of Proposition~\ref{prop:dual}.
  }
  \end{remark}

 \section{Alternative approaches and proofs}
\label{sec:alternative}

    An elegant and versatile analysis of the combinatorial structure of Lipschitz polytopes,  conducted by Gordon and Petrov  in \cite{gp}, can be with little care (but without introducing any new ideas) extended to the case of quasi-metrics.

   \medskip
   This fact, as kindly pointed by an anonymous referee,  provides a method for describing a large class of  quasi-metrics  which are combinatorially of ``cyclohedral type''.

     \medskip
   Here we give an outline of this method. (This whole section can be seen as a short  addendum to the paper \cite{gp}.)

  \bigskip
  A {\em combinatorial structure} on the (dual) pair of polytopes $Lip(\rho)$ and $KR(\rho)$ is, following \cite[Definition~2]{gp}, the collection of directed graphs $\mathcal{D}(\rho) = \{D(\alpha) \mid \alpha \mbox{ {\rm is a face of} } KR(\rho)\}$, where for each face $\alpha$ of $KR(\rho)$,
   \[
      (i,j) \in D(\alpha) \quad  \Leftrightarrow  \quad   \frac{e_j-e_i}{\rho(i,j)} \in \alpha.
   \]
   Following \cite[Definition 1]{gp}, a quasi-metric is {\em generic} if the triangle inequality is strict ($x\neq y\neq z  \Rightarrow \rho(x,z) < \rho(x,y) + \rho(y,z)$) and the polytope $KR(\rho)$ is simplicial ($Lip(\rho)$ is simple).

   \medskip
   In the case of a generic quasi-metric, the combinatorial structure   $\mathcal{D}(\rho)$ is a simplicial complex whose face poset is isomorphic to the face poset of the polytope $KR(\rho)$ (see Corollary~1 and Theorem~4 in \cite{gp}). Moreover, in this case $D(\alpha)$ is a directed forest (such that either the in-degree or the out-degree of each vertex is zero), and in particular if $\alpha$ is a facet then $D(\alpha)$ is a directed tree.

\medskip
Following \cite[Theorem 3]{gp} (see also  \cite[Theorem~4]{gp}) a directed tree (forest) $T$ is in $\mathcal{D}(\rho)$ if and only if it satisfies a {\em ``cyclic monotonicity''} condition (inequality (1) in \cite[Theorem 3]{gp}), indicating that $T$ can be thought of as an `optimal transference plan' for the transport of the corresponding measures.

\bigskip
It was shown in Section~\ref{sec:cyclo-prelim} (Proposition~\ref{prop:tree}) that the face poset of a cyclohedron can be also described as a poset of directed trees (corresponding to the diagrams of oriented arcs, as exemplified by Figures~\ref{fig-arc-2} and \ref{fig-arc-1}).

\medskip
From these observations arises a general plan for finding a generic quasi-metric $\rho$ such that the associated  combinatorial structure   $\mathcal{D}(\rho)$ is precisely the collection of trees associated to a cyclohedron. Moreover, this approach allows us (at least in principle) to characterize all generic quasi-metrics of ``combinatorial cyclohedral type''.

\medskip
Indeed, if  $\rho = (\rho(i,j))_{1\leq i,j \leq n}$ is an unknown quasi-metric (ranging over the space of all quasi-metric matrices), then one can characterize quasi-metrics of cyclohedral type by writing all inequalities of the type (1) in \cite[Theorem~3]{gp} (see also the simplification provided by \cite[Theorem 4]{gp}).

\begin{remark}{\rm
  Guided by the form of the metric $\rho$, described by the formula (\ref{eqn:explicit-d}) (Proposition~\ref{prop:quasi-metric-cyclo}),  the referee observed that the quasi-metric $\rho_\epsilon := d - \epsilon\cdot d^2$ (where $d$ is the clock quasi-metric and $\epsilon>0$ a sufficiently small number), is a good candidate for a cyclohedral quasi-metric.  (The details of related calculations will appear elsewhere.)
  }
\end{remark}

\begin{remark}{\rm
  The quasi-metric $\rho$ introduced in Proposition~\ref{prop:quasi-metric-cyclo} is somewhat exceptional since in this case we guarantee that $KR(\rho)$ is geometrically (affinely) and not only combinatorially (via face posets) isomorphic to a dual of a standard cyclohedron. This has some interesting consequences, for example this metric has the  property that the associated Lipschitz polytope $Lip(\rho)$ is a {\em Delzant polytope}. }
\end{remark}

\section{Concluding remarks}\label{sec:concluding}

\subsection{The result of Gordon and Petrov}
\label{sec:GP-concl}

The motivating result of J.~Gordon and F.~Petrov  \cite[Theorem~1]{gp} says that for a generic metric $\rho$ on a set of size $n+1$, the number of $(n-i)$-dimensional faces of the associated Lipschitz polytope (the dual of $KR(\rho)$, see the equation (\ref{eqn:Lip})) is equal to,
\[
   f_{n-i}(Lip(\rho)) = {{n+i}\choose{i,i,n-i}} = \frac{(n+i)!}{i!i!(n-i)!}.
\]
The link with the combinatorics of cyclohedra, established by Theorems~\ref{thm:glavna-A} and \ref{thm:glavna-B}, allows us to deduce this result from the known calculations of $f$-vectors of these polytopes. For example R.~Simion in \cite[Proposition~1]{si} proved that,
\[
     f_{i-1}(W_n) = {n \choose i}{n+i \choose  i} = \frac{(n+i)!}{i!i!(n-i)!}.
\]
Moreover, in light of Theorems~\ref{thm:glavna-A} and \ref{thm:glavna-B}, the generating series for these numbers have a new interpretation as a solution of a concrete partial differential equation, see for example
\cite[Sections 1.7. and 1.8.]{bp}.

\subsection{Tight hypergraphs}
\label{sec:tight}

The relationship between the cyclohedron and the (dual of the) root polytope is explained in Section~\ref{sec:revisited} as a special case of the relationship between tight hypergraphs $\mathcal{F}$ and their building closures $\mathcal{\widehat{F}}$.
For this reason it may be interesting to search for other examples of `tight pairs' $(\mathcal{\widehat{F}}, \mathcal{F})$ of hypergraphs.

\begin{exam}
  \label{prop:gen-ciklo-irred}{\rm
Let $C_n$ be the cycle graph on $n$ vertices (identified with their labels [n]) and let $\prec$ be the associated (counterclockwise) cyclic order on $[n]$.  For each ordered pair $(i, j)$ of indices let $[i,j]^0:= \{i, i+1, \dots, j\}$ be the associated `cyclic interval'.  For each $2\leq k\leq n-1$ let $BC_n^k\subset 2^{[n]}\setminus\{\emptyset\}$ be the hypergraph defined by,
\begin{equation}
BC_n^k := \{  [1,k]^0, [2,k+1]^0, \dots, [n, k-1]^0\}.
\end{equation}
Then $BC_n^k$ is a tight hypergraph on $[n]$. Moreover, if $\widehat{\mathcal{F}}_{n,k}: =\widehat{BC}_n^k$ is its building closure and $Q_{n,k}:=\Delta_{\widehat{\mathcal{F}}_{n,k}}$ the associated simple polytope than  for $k\neq k'$ the associated polytopes $Q_{n,k}$ and $Q_{n,k'}$ are combinatorially non-isomorphic.  }
\end{exam}

Note that for $k=2$ we recover the tight hypergraph described in Proposition~\ref{prop:ciklo-irred} (equation (\ref{eqn:1-2-3}))  and  in this case $Q_{n,2} = W_n$. Moreover observe that,
\[
     \mathcal{N}_{n,2} \varsupsetneq \mathcal{N}_{n,3} \supsetneq \dots \supsetneq   \mathcal{N}_{n,n-1},
\]
where $\mathcal{N}_{n,k}$ is the poset of nested sets in $\widehat{\mathcal{F}}_{n,k}$, which immediately implies that  $Q_{n,k} \ncong Q_{n,k'}$ for $k\neq k'$.

\subsection{Canonical quasitoric manifold over a cyclohedron}
\label{sec:canonical-quasitoric}

The cyclohedron $W_n$, together with the associated canonical map $\phi_n$, restricted to the set of vertices of $W_n^\circ$ ($\leftrightarrow$ the set of facets of $W_n$), defines a {\em combinatorial quasitoric pair} $(W_n, \phi_n)$ in the sense of  \cite[Definition~7.3.10]{bp}. Indeed, if $F_1,\dots, F_{n-1}$ are distinct facets of $W_n$ such that $\cap_{i=1}^{n-1}~F_i \neq\emptyset$, then the corresponding dual vertices $v_1, \dots, v_{n-1}$ of $W_n^\circ$ span a simplex and the vectors $\phi(v_1), \dots, \phi(v_{n-1})$ form a basis of the associated type A root lattice $\Lambda_n \cong \mathbb{Z}^{n-1}$ (spanned by the vertices of the root polytope $Root_n$).

\medskip
We refer to the associated quasitoric manifold $ M = M_{(W_n, \phi_n)}$ as the {\em canonical quasitoric manifold} over a cyclohedron $W_n$.

\subsection{The cyclohedron and the self-linking knot invariants}

It may be expected that the combinatorics of the map $\phi : W_n^\circ \rightarrow Root_n$, as illustrated by Theorems~\ref{thm:glavna-A} and \ref{thm:glavna-B} (and their proofs), may be of some relevance for other applications where the  cyclohedron $W_n$ played an important role. Perhaps the most interesting is the role of the cyclohedron in the combinatorics of the self-linking knot invariants (Bott and Taubes \cite{bt}, Voli\' c \cite{vo}). Other potentially interesting applications include  some problems of discrete geometry, as exemplified by  the `polygonal pegs problem' \cite{vz} and its relatives.

\bigskip\noindent
 {\bf Acknowledgements:}   The project was initiated during the program `Topology in Motion', \url{https://icerm.brown.edu/programs/sp-f16/}, at the {\em Institute for Computational and Experimental Research in Mathematics} (ICERM, Brown University). With great pleasure we acknowledge the support, hospitality and excellent working conditions at ICERM. The research of Filip Jevti\'{c} is a part of his PhD project at the University of Texas at Dallas, performed under the supervision and with the support of Vladimir Dragovi\'{c}.  We would also like to thank the referee for very useful comments and suggestions, in particular  for the observations incorporated in Section~\ref{sec:alternative}.

\end{document}